\documentclass[12pt,leqno]{amsart}
\usepackage{a4wide}
\usepackage{amssymb}
\usepackage{amsmath}
\usepackage{amsthm}
\usepackage{verbatim}
\usepackage{xypic}
\numberwithin{equation}{section}

\theoremstyle{plain}
\newtheorem{theo}{Theorem}[section]
\newtheorem{lem}[theo]{Lemma}
\newtheorem{prop}[theo]{Proposition}
\newtheorem{cor}[theo]{Corollary}

\newtheorem{conj}[theo]{Conjecture}
\theoremstyle{remark}
\newtheorem{rem}[theo]{Remark}

\theoremstyle{definition}

\newcommand{\ff}{\hbox{if }}
\newcommand{\Sym}{\hbox{Sym}}

\renewcommand{\Im}{\operatorname{Im}}
\newcommand{\A}{{\mathbb A}}
\newcommand{\C}{{\mathbb C}}

\newcommand{\Q}{{\mathbb Q}}
\newcommand{\Z}{{\mathbb Z}}

\newcommand{\zxz}[4]{\begin{pmatrix} #1 & #2 \\ #3 & #4 \end{pmatrix}}

\newcommand{\kzxz}[4]{\left(\begin{smallmatrix} #1 & #2 \\ #3 & #4\end{smallmatrix}\right) }
\newcommand{\kabcd}{\kzxz{a}{b}{c}{d}}

\newcommand{\Gal}{\operatorname{Gal}}

\newcommand{\CL}{\mathcal C\mathcal L}

\renewcommand{\div}{\operatorname{div}}

\newcommand{\CM}{\mathcal{CM}}

\newcommand{\norm}{\operatorname{N}}
\newcommand{\End}{\operatorname{End}}
\newcommand{\Hom}{\operatorname{Hom}}

\newcommand{\disc}{\operatorname{disc}}

\newcommand{\Ht}{\hbox{ht}}

\renewcommand{\O}{\mathcal O}
\newcommand{\GL}{\operatorname{GL}}

\newcommand{\tr}{\operatorname{tr}}
\newcommand{\ord}{\operatorname{ord}}
\newcommand{\SL}{\operatorname{SL}}

\newcommand{\cha}{{\text{\rm char}}}
\newcommand{\diag}{{\text{\rm diag}}}
\newcommand{\Ind}{\operatorname{Ind}}
\newcommand{\Aut}{\operatorname{Aut}}
\newcommand{\CalZ}{\mathcal Z}
\newcommand{\CalM}{\mathcal M}
\newcommand{\Pic}{\operatorname{Pic}}
\newcommand{\CH}{\operatorname{CH}}
\newcommand{\Fal}{\operatorname{Fal}}
\newcommand{\Pet}{\operatorname{Pet}}
\newcommand{\pr}{\operatorname{pr}}
\newcommand{\Cm}{\operatorname{CM}}

\begin{document}

\title{An arithmetic intersection formula on Hilbert modular surfaces}
\author[ Tonghai Yang]{Tonghai Yang}
\address{Department of Mathematics, University of Wisconsin Madison, Van Vleck Hall, Madison, WI 53706, USA}
\email{thyang@math.wisc.edu} 
\thanks{ Partially supported by NSF grants DMS-0302043, 0354353,
 and a Chinese NSF grant NSFC-10628103.}

\begin{abstract}
  In this paper, we obtain an explicit arithmetic intersection formula on a Hilbert modular surface between the diagonal
  embedding of the modular curve  and a CM cycle associated to a
  non-biquadratic CM quartic field. This confirms a special case of
  the author's conjecture with J. Bruinier in \cite{BY}, and is a
  generalization of the beautiful factorization formula of
  Gross and Zagier on singular moduli.  As an application, we proved
  the first non-trivial non-abelian Chowla-Selberg formula, a
  special case of Colmez conjecture. 
\end{abstract}

  \maketitle

\section{Introduction}

Intersection theory and Arakelov theory  play important roles
in algebraic geometry 
and  number theory. Indeed, some of the deepest results and conjectures,
such as
Faltings's proof of
Mordell's Conjecture, and the work of Gross and Zagier on the Birch and
Swinnerton-Dyer Conjecture, highlight these roles.
Deep information typically follows from the
derivation of explicit intersection formulae.
For example, consider the
Gross-Zagier formula \cite{GZ} and  its generalization
by Shou-Wu Zhang
\cite{Zh1}, \cite{Zh2}, \cite{Zh3}.
We also have
recent work on an arithmetic
Siegel-Weil formula by Kudla, Rapoport, and the author
\cite{KuAnnal}, \cite{KRY1}, \cite{KRY2}, along with
work of Bruinier, Burgos-Gil, and K\"uhn  on
an arithmetic Hirzebruch-Zagier formula \cite{BBK}.
There are many other famous examples of explicit intersection
formulae.
There is the work of
Gross and Zagier on singular moduli \cite{GZ2},
the work of Gross and Keating on modular polynomials \cite{GK},
along with its many applications
(for example, see  \cite{KuAnnal}, \cite{KR1}, \cite{KR2}),
as well as the recent results of Kudla and Rapoport
\cite{KR1, KR2} in the context of
Hilbert modular surfaces and Siegel modular
$3$-folds.

In all of these works, the intersecting cycles are symmetric and
are of similar type. We investigate two {\it different types} of
cycles in a Hilbert modular surface  defined over $\mathbb Z$:
arithmetic Hirzebruch-Zagier divisors, and arithmetic CM cycles
associated to  non-biquadratic quartic CM fields. These cycles
intersect properly, and in earlier work with Bruinier \cite{BY},
we conjectured the corresponding  arithmetic intersection formula.
The truth of this formula has applications  to a well known
conjecture of Colmez which aims to generalize the classical
Chowla-Selberg formula \cite{Col}, as well as a conjecture of
Lauter on the denominators of the evaluations of Igusa invariants
at CM points \cite{La}. Here we prove a special case of the
conjectured formula, and as a consequence we obtain the first
generalization of the Chowla-Selberg formula to {\it non-abelian
CM number fields}. This result  confirms Colmez's conjecture in
this case. It also confirms Lauter's conjecture in certain cases,
but for brevity we shall omit a detailed discussion.

We begin by fixing notation.
Let $D\equiv 1 \mod 4$ be  prime, and let $F=\mathbb
Q(\sqrt D)$ with the ring of integers $\O_F=\mathbb
Z[\frac{D+\sqrt D}2]$ and different $\partial_F =\sqrt D \O_F$.
Let  $\mathcal M$ be  the Hilbert moduli stack over $\mathbb Z$
representing the moduli problem that  assigns a base scheme  $S$
over $\mathbb Z$ to the set of the triples $(A, \iota, \lambda)$,
where (\cite[Chapter 3]{Go} and \cite[Section 3]{Vo})

(1)\quad  $A$ is a abelian surface over $S$.

(2) \quad $\iota: \O_F \hookrightarrow \End_S(A)$ is real
multiplication of $\O_F$ on $A$.

(3) \quad $\lambda: \partial_F^{-1} \rightarrow
P(A)=\Hom_{\O_F}(A, A^\vee)^{\hbox{sym}}$ is a
$\partial_F^{-1}$-polarization (in the sense of Deligne-Papas)
satisfying the condition:
\begin{equation} \label{eq1.1}
 \partial_F^{-1}\otimes A  \rightarrow A^\vee, \quad  r \otimes a
 \mapsto \lambda(r)(a)
 \end{equation}
 is an isomorphism (of Abelian schemes).

 Next, for an integer $m \ge 1$, let $\mathcal T_m$ be the integral Hirzebruch-Zagier divisors in $\mathcal M$
 defined in \cite[Section 5]{BBK}, which is the flat
 closure of the classical Hirzebruch-Zagier divisor $T_m$ in $\mathcal
 M$. For $m=1$, $\mathcal T_1$ has the following simple moduli
 description.
Let $\mathcal E$ be the moduli stack over $\mathbb Z$ of elliptic
curves, then $E \mapsto (E\otimes \O_F, \iota, \lambda)$ is a closed
immersion from $\mathcal E$ into $\mathcal M$, and its  image is
$\mathcal T_1$.
$$
\iota: \O_F  \hookrightarrow \End_S(E) \otimes \O_F =\End_{S \otimes
\O_F}(E\otimes \O_F) \hookrightarrow \End_S(E \otimes \O_F)
$$
is the natural embedding, and
$$
\lambda: \partial_F^{-1} \rightarrow \Hom_{S \otimes \O_F} (E
\otimes \O_F, E \otimes \partial_F^{-1})^{\hbox{sym}}, \quad
\lambda(z) (e\otimes x) = e \otimes xz.
$$
By abuse of notation, we will identify $\mathcal E$ with $\mathcal
T_1$.

Finally,  let $K=F(\sqrt\Delta)$ be a quartic non-biquadratic CM
number field with real quadratic subfield $F$.  Let $\CM(K)$ be
the moduli stack over $\mathbb Z$ representing the moduli problem
which  assigns a base scheme $S$ to the set of  the triples $(A,
\iota, \lambda)$ where $\iota: \O_K \hookrightarrow \End_S(A)$ is
an CM action of $\O_K$ on $A$, and $ (A, \iota|_{\O_F},
\lambda)\in \mathcal M(S)$ such that the Rosati involution
associated to $\lambda$ induces to the complex conjugation of
$\O_K$. The map $(A, \iota, \lambda) \mapsto (A, \iota|_{\O_F},
\lambda)$ is a finite proper map  from $\CM(K)$ into $\mathcal M$,
and we denote its direct image in $\mathcal M$ still by $\CM(K)$
by abuse of notation. Since $K$ is non-biquadratic, $\mathcal T_m$
and $\CM(K)$ intersect properly. A basic question is to compute
their arithmetic intersection number (see Section \ref{newsect2}
for definition).  We have the following conjectured intersection
formula, first stated in \cite{BY}. To state the conjecture, let
$\Phi$ be a CM type of $K$ and let $\tilde K$ be reflex field of
$(K, \Phi)$. It is also a quartic non-biquadratic CM field with
real quadratic field $\tilde F= \mathbb Q(\sqrt{\tilde D})$ with
$\tilde D =\Delta \Delta'$. Here $\Delta'$ is the Galois conjugate
of $\Delta$ in $F$.

\begin{conj} (Bruinier and Yang \cite{BY}) \label{conj}  Let the notation be as above.
Then
\begin{equation}
\mathcal T_m .\CM(K) = \frac{ 1}{2} b_m
\end{equation}
or equivalently
\begin{equation}\label{conjl}
(\mathcal T_m .\CM(K))_p = \frac{1}{2} b_m(p)
\end{equation}
for every prime $p$. Here $$ b_m =\sum_p b_m(p) \log p$$ is defined
as follows:
\begin{equation} \label{eqI1.3}
b_m(p) \log p =\sum_{\mathfrak p|p} \sum_{t=\frac{n+m\sqrt{\tilde
D}}{2D} \in d_{\tilde K/\tilde F}^{-1}, |n| <m \sqrt{\tilde D}}
B_t(\mathfrak p)
\end{equation}
where \begin{equation} B_t(\mathfrak p)
 =\begin{cases}
  0 &\ff \mathfrak p \hbox{ is split in} \tilde K,
  \\
  (\ord_{\mathfrak p} t_n+1)
\rho(t d_{\tilde K/\tilde F} \mathfrak p^{-1}) \log |\mathfrak p|
 &\ff \mathfrak p \hbox{ is not split in} \tilde K,
\end{cases}
\end{equation}
and
$$
\rho(\mathfrak a) =\#\{ \mathfrak A \subset \O_{\tilde K} :
N_{\tilde K/\tilde F} \mathfrak A =\mathfrak a\}.
$$
\end{conj}

Notice that the conjecture implies that $\mathcal T_m.\CM(K) =0$
unless $4Dp| m^2 \tilde D -n^2$ for some integer $0 \le  n  < m\sqrt
{\tilde D}$, in particular one has to have $p \le \frac{m^2 \tilde
D}{4 D}$.

Throughout  this paper, we assume that $ K$ satisfies the following
condition---we call it condition $(\clubsuit)$:

\begin{equation} \label{eqOK}
\O_K =\O_F + \O_F \frac{w+\sqrt\Delta}2
\end{equation}
is free over $\O_F$ and that  $\tilde D =\Delta \Delta' \equiv 1
\mod 4$ is square free ($w \in \O_F$). Under this assumption, one can
show that $d_K= D^2 \tilde D$, and $d_{\tilde K} =\tilde D^2 D$,
and $N d_{\tilde K/\tilde F} =D$. Here $d_K$ is the discriminant
of $K$, and  $d_{\tilde K/\tilde F}$ is the relative discriminant
of $\tilde K/\tilde F$. The main purpose of this paper is to prove
the conjecture when $m=1$, and to give a simple procedure for computing
$b_1(p)$.

\begin{theo} \label{maintheo} Under the condition $(\clubsuit)$,
Conjecture \ref{conj} holds for $m=1$. 
\end{theo}

 We prove the theorem by computing the local intersection $(\CM(K).\mathcal T)_p$
and $b_1(p)$ at given $p$ separately and comparing them. On the
geometric side, to
 a  geometric intersection point $\iota: \O_K \hookrightarrow
\End(E)\otimes \O_F$ we first associate   a positive integer $n$,
a sign $\mu=\pm 1$, and a $2\times 2$ integral matrix $T(\mu n)$
with $\det T(\mu m) = \frac{\tilde D -n^2}{D} \in 4 p \mathbb
Z_{>0}$ (Proposition \ref{propI2.2}). Next, we use Gross and
Keating's beautiful formula \cite{GK} to show the local
intersection index at the  geometric point $\iota$ is equal to
$\frac{1}{2}(\ord_{p} \frac{\tilde D-n^2}{4D} +1)$, depending only
on $T(\mu n)$, not on the geometric point itself (Theorem
\ref{propI0.4}). Practically, the local intersection index $\iota$
is  the largest integer $m$ this action can be lifted to $W/p^m$
where $W$ is the Witt ring of $\bar{\mathbb F}_p$. The
independence on the geometric points is essential and leads us to
a simpler problem of counting the number of geometric points
$\iota: \O_K \hookrightarrow \End(E)\otimes \O_F$ whose associated
matrices is $T(\mu n)$, which is a   local density problem
representing $T(\mu n)$ by a ternary integral lattice. Explicit
computation for the local density problem is given in
\cite{YaDensity1} and \cite{YaDensity2}, but the formula at $p=2$
is extremely complicated in general. We circumvent it in this
special case by switching it to similar local density problem with
clean known answer in Section \ref{sect3}, and obtain the
following intersection formula.

\begin{theo} \label{newtheo1.3}  Let the notation and assumption be as in  Theorem
\ref{maintheo}, and let $p$ be a prime number. Then
 \begin{equation} \label{eq1.7}
(\mathcal T_1.\CM(K))_p = \frac{1}{2} \sum_{0 < n <\sqrt{\tilde
D}, \, \frac{\tilde D -n^2}{4D} \in p
 \mathbb Z_{>0}}(\ord_p \frac{\tilde D-n^2}{4D} +1) \sum_{\mu }\beta(p, \mu n),
\end{equation}
where
$$ \beta(p, \mu n) = \prod_{l | \frac{\tilde D
-n^2}{4D}} \beta_l(p,\mu n)
$$
is given as follows. Given a positive integer $0<n <\sqrt{\tilde D}$
with $\frac{\tilde D-n^2}{4D} \in p\mathbb Z_{>0}$ as in
(\ref{eq1.7}), there is one sign  $\mu =\pm 1$ (both signs if $D|n$)
and a unique positive definite integral $2\times 2$ matrix $T(\mu
n)$ satisfying the conditions in Lemma \ref{lemI2.1}. For a fixed
prime $l$, $T(\mu n)$ is $\GL_2(\mathbb Z_l)$-equivalent to
$\diag(\alpha_l, \alpha_l^{-1} \det T(\mu n))$ with $\alpha_l \in
\mathbb Z_l^*$. Let $t_l =\ord_l \frac{\tilde D -n^2}{4D}=\ord_l
T(\mu n)- 2 \ord_l 2$. Then
$$
\beta_l(p, \mu n) =\begin{cases}
   \frac{1- (-\alpha_p, p)_p^{t_p}}2 &\ff l=p,
   \\
   \frac{ 1+ (-1)^{t_l}}2 &\ff l \ne p, (-\alpha_l, l)_l =-1,
  \\
   t_l + 1 &\ff   l \ne p, (-\alpha_l, l)_l =1.
\end{cases}
$$
\end{theo}

The theorem has the following interesting consequence.

\begin{cor} \label{cor1.3} Assume $(\clubsuit)$ and $\tilde D < 8D$. Then
$\mathcal T_1.\CM(K)=0$, i.e., there is no elliptic curve $E$ such
that $E \otimes \O_F$ has CM by $\O_K$.
\end{cor}

 In Section \ref{sect4},  we compute $b_1(p)$
and show that it equals twice the right hand side of (\ref{eq1.7})
and thus prove Theorem \ref{maintheo}. From the definition, it is
sufficient to prove an identity for each positive integer $n$ with
$\frac{\tilde D -n^2}{4D} \in p\mathbb Z_{>0}$. After some
preparation, one sees that the key is to relate whether $\tilde
K/\tilde  F$ is split or inert at a prime $\mathfrak l$ to the
local property of $T(\mu n)$ at prime $l= \mathfrak l \cap \mathbb
Z$.  We prove this unexpected connection in Lemma \ref{lemI4.2},
and finish the computation of $b_1(p)$  in Theorem \ref{theoI4.4}.

  It is worth noting a
mysterious identity underlining the conjecture. On the one hand,
it is clear from our proof and a general program of Kudla
\cite{Ku2} that the intersection number is summation over some
Fourier coefficients of the central derivative of some incoherent
Siegel-Eisenstein series of genus $3$. On the other hand, it is
clear from \cite{BY} that $b_m(p)$ comes from summation of certain
Fourier coefficients of the central derivative of incoherent
Eisenstein series on Hilbert modular surface. Viewing this
identity as an identity relating the two seemingly unrelated
Eisenstein series, one can naturally ask whether it is a pure
accident, or there is some hidden gem?

 Now we briefly describe an application of   Theorem
 \ref{maintheo} to a conjecture of Colmez,
 which is a beautiful generalization of the celebrated Chowla-Selberg formula.
  In
 proving the famous Mordell conjecture, Faltings introduces the
 so-called Faltings height $h_{\hbox{Fal}}(A)$ of an Abelian variety
 $A$, measuring the complexity of $A$ as a point in  a Siegel modular
 variety. When $A$ has complex multiplication, it  only depends on
 the CM type of $A$ and  has a simple description as follows. Assume that $A$ is defined over
 a number field $L$ with good reduction everywhere, and
 let $\omega_A \in \Lambda^g \Omega_A$ be a Neron differential
 of $A$ over $\O_L$, non-vanishing everywhere, Then
 the Faltings' height of $A$ is defined as (our normalization is
 slightly different from that of \cite{Col})
 \begin{equation}
 h_{\hbox{Fal}}(A) =-\frac{1}{2[L:\Q]}
    \sum_{\sigma: L \hookrightarrow \mathbb C} \log
     \left|(\frac{1}{2 \pi i})^g \int_{ \sigma(A)(\C)} \sigma(\omega_A) \wedge
     \overline{\sigma(\omega_A)}\right| + \log \# \Lambda^g \Omega_A/\O_L\omega_A.
     \end{equation}
Here $g=\dim A$.  Colmez gives a beautiful conjectural formula to
 compute the Faltings height of a CM abelian variety
 in terms of the log derivative of certain Artin L-series associated to the CM type  \cite{Col},
 which is consequence of   his  product formula conjecture of $p$-adic
 periods in the same paper. When $A$ is a CM elliptic curve, the
 height conjecture is a reformulation of  the well-known Chowla-Selberg
 formula relating the CM values of the usual  Delta function
 $\Delta$ with the values  of  the  Gamma function  at rational
numbers. Colmez proved his conjecture  up to a multiple of $\log 2$
when the CM field (which acts on $A$) is abelian, refining Gross's
 \cite{Gr} and Anderson's  \cite{Ad} work. A key point is that such CM abelian varieties are
quotients of the Jacobians of the Fermat curves, so one has a model
to work with. When the CM number field is non-abelian,  nothing is
known. Conjecture \ref{conj}, together with \cite[Theorem 1.4]{BY},
would prove Colmez's conjecture for non-biquadratic quartic CM
fields, confirming the first {\it non-abelian} case. More precisely,
let $K$ be a non-biquadratic CM number field with totally real
quadratic subfield $F=\mathbb Q(\sqrt D)$. Let $\chi$ be the
quadratic Hecke character of $F$ associated to $K/F$ by the global
class field theory, and let
\begin{equation}
\Lambda(s, \chi) = C(\chi)^{\frac{s}2}
\pi^{-s-1}\Gamma(\frac{s+1}2)^2 L(s, \chi)
\end{equation}
be the complete L-function of $\chi$ with $C(\chi) =D N_{F/\mathbb
Q} d_{K/F}$. Let
\begin{equation}
\beta(K/F)
 = \frac{\Gamma'(1)}{\Gamma(1)}
-\frac{\Lambda'(0, \chi )}{\Lambda(0, \chi )} -\log4\pi .
\end{equation}
In this case, the conjectured formula of Colmez on the Faltings's
height of a CM abelian variety $A$ of type $(K, \Phi)$ does not even
depend on the CM type $\Phi$ and is given by (see \cite{Ya3})
\begin{equation} \label{eqCol}
h_{\hbox{Fal}}(A)=\frac{1}2 \beta( K/ F) .
\end{equation}

In Section \ref{newsect7}, we will prove using Theorem
\ref{maintheo}, and \cite[Theorem 1.4]{BY}.

\begin{theo} \label{theo1.4} Let $K$ be a non-biquadratic CM quartic CM field of
discriminant $D^2 \tilde D$ with $D =5, 13, \hbox{ or } 17$, and
$\tilde D \equiv 1 \mod 4$ prime. Then Colmez's conjecture
(\ref{eqCol}) holds.
\end{theo}

Theorem \ref{maintheo} also has implications for
Lauter's conjecture
on the denominator of Igusa invariants at CM points and bad
reduction of CM genus two curves in the special cases $D=5, 13$, and
$17$. To keep this paper short, concise, and to the point, we omit
this application and refer the reader to \cite{Ya4} for this
application, where we prove Conjecture \ref{conj} under the
condition (\ref{eqOK}) and that $\tilde D \equiv 1 \mod 4$ is a
prime. The idea is to prove a weaker version  of the conjecture for
$\mathcal T_q$ when $q$ is a prime split in $F=\mathbb Q(\sqrt D)$
(up to a multiple of $\log q$), and then combining it with
\cite[Theorem 1.4]{BY} and \cite[Theorem 4.15]{BBK} to derive the
general case. Although the proof of the weaker version is similar to
the case $m=1$ in this paper in principle, the argument is much more
complicated and needs new ideas. The first difficulty is that
instead of simple $\End_{\O_F}( E\otimes \O_F)= \End(E) \otimes
O_F$, $\End_{\O_F}(A)$ does not have a good global interpretation.
So we have to work locally in terms of Tate modules and Dieudonne
modules.  Second, the local density problem is no longer a problem
representing one matrix by a lattice. Instead, it is really a local
Whittaker integral. We have to use a totally different method to
compute the local integral.

  Here is the organization of this paper. In Section
  \ref{newsect2}, we give basic definition for arithmetic
  intersection  and Faltings' heights in stacks, following
  \cite{KRY2}. We also show that $\mathcal T_1$ is isomorphic to
  the stack of elliptic curves. In Section \ref{newsect3}, we
  briefly sketch a proof of Theorem \ref{maintheo} in the
  degenerate case $D=1$, which is also a new proof of the
  Gross-Zagier formula on factorization of singular moduli
  \cite{GZ2}.  In Section \ref{sect2}, we  use a beautiful formula
  of Gross and Keating \cite{GK} to compute the local intersection
  index of $\mathcal T_1$ and $\CM(K)$ at a geometric intersection point. In
  Section \ref{sect3}, we count the number of geometric intersection
  points of $\mathcal T_1$ and $\CM(K)$ and prove Theorem
  \ref{newtheo1.3}. In Section \ref{sect4}, we compute $b_1(p)$
  and finish the proof of Theorem \ref{maintheo}. In the last
  section, we prove Theorem \ref{theo1.4}.

{\bf Acknowledgement:} The author thanks Bruinier, Kudla, K\"uhn,
Lauter,  Olsson, Ono,  Rapoport,  Ribet, and Shou-Wu Zhang  for
their help during the preparation of this paper. He thanks the
referee for his/her careful reading of this paper and very helpful
suggestions which makes the exposition of this paper much better.
Part of the work was done when the author visited the Max-Planck
Instit\"ut of Mathematik at Bonn, MSRI, the  AMSS and the
Morningside Center of Mathematics  at Beijing. He thanks these
institutes for providing him wonderful working environment.

\section{Basic definitions} \label{newsect2}

We basically follow \cite[Chapter 2]{KRY2} in our definition of
arithmetic intersection and Faltings' height on DM-stacks which
have a quotient presentation.

Let $\mathcal M$ be a regular DM-stack of dimension $n$ which is
proper and flat over $\mathbb Z$. Two cycles $\mathcal Z_1$ and
$\mathcal Z_2$ in $\mathcal M$ of co-dimensions $p$ and $q$
 respectively with $p+q=n$ intersect properly if $\mathcal Z_1
 \cap \mathcal Z_2 =\mathcal Z_1 \times_{\mathcal M} \mathcal Z_2$
 is a DM-stack of dimension $0$. In this case, we define the
 (arithmetic) intersection number as
 \begin{equation}
 \mathcal Z_1.\mathcal Z_2=\sum_{p} \sum_{x \in \mathcal Z_1
 \cap \mathcal Z_2(\bar{\mathbb F}_p) }\frac{1}{\#\Aut(x)} \log \#
 \tilde{\O}_{\CalZ_1\cap \CalZ_2, x}
    =\sum_{p} \sum_{x \in \mathcal Z_1
 \cap \mathcal Z_2(\bar{\mathbb F}_p) } \frac{1}{\#\Aut(x)}
 i_p(\CalZ_1, \CalZ_2, x) \log p
 \end{equation}
 where $\tilde{\O}_{\CalZ_1\cap\CalZ_2, x}$ is the strictly local
 henselian ring of $\CalZ_1\cap\CalZ_2$ at $x$,
 $$
 i_p(\CalZ_1, \CalZ_2, x) = \operatorname{Length} \tilde{\O}_{\CalZ_1\cap\CalZ_2, x}
 $$
 is the local intersection index of $\CalZ_1$ and $\CalZ_2$ at $x$. If $\phi: \CalZ \rightarrow \CalM$ is a finite proper and
 flat map from stack $\CalZ$ to $\CalM$, we will identify  $\CalZ$
 with its direct image $\phi_*\CalZ$ as a cycle of $\CalM$, by abuse
 of notation.

  Now we further assume that its generic fiber $M=\CalM_{\mathbb C}
  =[\Gamma\backslash X]$ is a  quotient stack of a regular proper scheme $X$,  where $\Gamma$ is a finite
  group acting on $X$. Let
  $$
  \pr: X \rightarrow M
  $$
  be the natural projection. We define the arithmetic Picard group
  $\widehat{\Pic}(\CalM)$ and the arithmetic Chow group $\widehat{\CH}^1(\CalM)$
  as in  \cite[Chapter 2]{KRY2}.  For example, let
  $\hat{Z}^1(\CalM)$
  is $\mathbb R$-vector space generated by
  $(\CalZ, g)$, where $\CalZ$ is a prime divisor in $\CalM$ (a closed irreducible reduced substack of codimension $1$
  in $\CalM$ which is locally in \'etale topology by a Cartier divisor),
  and $g$ is a Green function for $Z=\CalZ(\mathbb C)$. It
  means the following.
   Let $\tilde Z=\pr^{-1}(Z)$ be the associated
  divisor in $X$. Then the Dirac current $\delta_Z$ on $M$ is given
  by
  $$
  < \delta_Z, f>_M =\frac{1}{\# \Gamma} <\delta_{\tilde Z}, f>_X
  $$
  for every $C^\infty$ function on $M$ with compact support (i.e.,
  every $\Gamma$-invariant $C^\infty$ function  on $X$ with compact
  support). A Green function for $Z$ is defined to be a $\Gamma$-invariant
  function $g$ for $\tilde Z$. In such a case, we also have
  naturally
  $$
  d d^c g + \delta_Z = [\omega]
  $$
  as currents in $M$ for some smooth $(1, 1)$-form $\omega$ on
  $M$---a $\Gamma$-invariant smooth $(1, 1)$-form on $X$
  (see \cite[(2.3.11)]{KRY2}). Although $n=1$ is assumed in \cite{KRY2},
  the same argument  holds for all $n$. For a rational function $f
  \in Q(\CalM)^*$, one defines
  $$
  \widehat{\div} (f) = (\div f, -\log |f|^2) \in \hat{Z}^1(\CalM)
  $$
Then $\widehat{\CH}^1(\CalM)$ is the quotient space of
$\hat{Z}^1(\CalM)$ by the $\mathbb R$-vector space generated by
$\widehat{\div}(f)$.

   There is a natural
  isomorphism
  $$
  \widehat{\Pic}(\CalM) \cong \widehat{\CH}^1(\CalM), \quad
  $$
  which is induced by $\hat{\mathcal L}=(\mathcal L, \|\, \|) \mapsto (\div s, -\log \|
  s\|^2)$, $s$ is a rational section of $\mathcal L$. Given a
  finite proper and flat map $\phi: \CalZ \rightarrow \CalM$, it
  induces a pull-back map
$\phi^*$ from $  \widehat{\CH}^1(\CalM)$ to
$\widehat{\CH}^1(\CalZ)$, and from $\widehat{\Pic}(\CalM)$ to
$\widehat{\Pic}(\CalZ)$.  When $\CalZ$ is a prime cycle of dimension
$1$ (codimension $n-1$), and $\hat{\mathcal L}$ is a metrized line
bundle on $\CalM$, we define the Faltings height
\begin{equation}
h_{\hat{\mathcal L}}(\CalZ) = \widehat{deg} (\phi^*\hat{\mathcal L})
\end{equation}
where $\phi$ is the natural embedding of $\CalZ$ to $\CalM$.  Here
the arithmetic degree on $\widehat{\Pic}(\CalZ)$ is defined as in
\cite[(2.18) and (2.19)]{KRY2}.  In particular, if $s$ is a
(rational) section of $\mathcal L$ such that $\div s$ intersects
properly with $\mathcal Z$, we have
\begin{equation}
h_{\hat{\mathcal L}}(\CalZ) = \CalZ.\div s - \sum_{z \in Z}
\frac{1}{\# \Aut(z)} \log ||s(z)||
\end{equation}
where $Z=\CalZ(\mathbb C)$. Equivalently, in terms of arithmetic
divisors,
\begin{equation}
h_{(\CalZ_1, g_1)}(\CalZ) = \CalZ_1.\CalZ + \frac{1}2 \sum_{z \in
Z} \frac{1}{\# \Aut(z)} g_1(z), \quad (\CalZ_1, g_1) \in
\widehat{\CH}^1(\CalM)
\end{equation}
if $\CalZ_1$ and $\CalZ$ intersect properly. The Faltings height
is a bilinear map on $\widehat{\CH}^1(\CalM) \times
Z^{n-1}(\CalM)$, which does not factor through $\CH^{n-1}(\CalM)$.

Now come back to our specific case. Let $F=\mathbb Q(\sqrt D)$ be
a real quadratic field with $D\equiv 1 \mod 4$ being prime.   Let
$\CalM$ be the Hilbert modular stack over $\mathbb Z$ defined in
the introduction. It is regular and  flat over $\mathbb Z$ but not
proper (\cite{DP}). Let $\tilde{\CalM}$ be a fixed Toroidal
compactification of $\CalM$, then $\tilde{\CalM}_{\mathbb C}$ and
$\CalM_{\mathbb C}$ have quotient presentation (e.g.,
$\CalM(\mathbb C) = [\Gamma \backslash Y(N)]$ with $Y(N)=
\Gamma(N) \backslash \mathbb H^2$, and $\Gamma
=\Gamma(N)\backslash \SL_2(\O_F)$ for $N \ge 3$). Let $K=
F(\sqrt\Delta)$ be a non-biquadratic quartic CM number field with
real quadratic  subfield $F$, and let $\CM(K)$ be the CM cycle
defined  in the introduction. Notice that $\CM(K)$ is closed in
$\tilde{\CalM}$. $K$ has four different CM types $\Phi_1$,
$\Phi_2$, $\rho \Phi_1=\{ \rho\sigma: \, \sigma \in \Phi_1\}$, and
$\rho\Phi_2$, where $\rho$ is the complex conjugation in $\mathbb
C$. If $x=(A, \iota, \lambda) \in \CM(K)(\mathbb C)$, then $(A,
\iota, \lambda)$ is a CM abelian surface over $\mathbb C$ of
exactly one CM type $\Phi_i$ in $\mathcal M(\mathbb C)=\SL_2(\O_F)
\backslash \mathbb H^2$ as defined in \cite[Section 3]{BY}. Let
$\Cm(K, \Phi_i)$ be set of (isomorphism classes) of CM abelian
surfaces of CM type $(K, \Phi_i)$ as in \cite{BY}, viewed as a
cycle in $\mathcal M(\mathbb C)$. Then it was proved in \cite{BY}
$$
\Cm(K) =\Cm(K, \Phi_1) + \Cm(K, \Phi_2) =  \Cm(K, \rho \Phi_1) +
\Cm(K, \rho \Phi_2)
$$
is defined over $\mathbb Q$. So we have

\begin{lem} \label{newlem2.1} One has
$$
\CM(K)(\mathbb C) = 2\Cm(K)
$$
in   $\CalM(\mathbb C)$.
\end{lem}

Next, recall that the Hirzebruch-Zagier divisor $T_m$  is given by
\cite{HZ}
$$
T_m(\mathbb C)=\SL_2(\O_F)\backslash \{ (z_1, z_2) \in \mathbb
H^2: (z_2,
1)A  \left(\substack{ z_1 \\
1}\right) =0 \hbox{ for  some } A \in L_m \},
$$
where
$$
L_m=\{ A =\kzxz {a} {\lambda} {\lambda'} {b} :\, a, b \in \mathbb
Z, \lambda \in
\partial_F^{-1}, ab -\lambda \lambda' =\frac{m}D\}.
$$
 $T_m$
is empty if $(\frac{D}m)=-1$. Otherwise, it is  a finite union of
irreducible curves and is actually defined over $\mathbb Q$.  In
particular, $T_1(\mathbb C)$ is the diagonal image of modular
curve $\SL_2(\mathbb Z) \backslash \mathbb H$ in $\mathcal
M(\mathbb C)$. Following \cite{BBK}, let $\mathcal T_m$ be the
flat closure of $T_m$ in $\mathcal M$.

\begin{lem} \label{newlem2.3} Let $\mathcal E$ be the moduli stack over $\mathbb Z$  of elliptic
curves.  Let $\phi: \mathcal E \rightarrow \mathcal M$ be given by
$\phi(E) =(E\otimes \O_F, \iota_F, \lambda_F)$ for any elliptic
curve over a base scheme $S$, where
$$
\iota_F: \O_F  \hookrightarrow \End_S(E) \otimes \O_F =\End_{\O_S
\otimes \O_F}(A) \subset \End_S(A)
$$
is the natural embedding, and
$$
\lambda_F: \partial_F^{-1} \rightarrow \Hom_{\O_E} (E \otimes
\O_F, E \otimes \partial_F^{-1})^{\hbox{sym}}, \quad  \lambda_F(z)
(e\otimes x) = e \otimes xz.
$$
Then $\phi$ is a closed immersion and $\phi(\mathcal E) = \mathcal
T_1$.
\end{lem}
\begin{proof} It is known \cite[Proposition 5.14]{BBK} that $\phi$ is a proper map and
its image is $\mathcal T_1$. To show it is a closed immersion as
stacks, it is enough to show
$$
\hbox{Isom}(E, E') \cong \hbox{Isom}(\phi(E), \phi(E')),\quad f
\mapsto \phi(f).
$$
Clearly, if $f: E \rightarrow E'$ is an isomorphism, $\phi(f)$ is an
isomorphism between $\phi(E)$ and $\phi(E')$. On the other hand, if
$g: \phi(E) \rightarrow \phi(E')$ is an isomorphism, i.e., $g:
E\otimes \O_F \rightarrow E'\otimes \O_F$ is an $\O_F$-isomorphism
such that
\begin{equation}
  g^{\vee} \circ \lambda_F(r) \circ g = \lambda_F(r)
  \end{equation}
  for any $r \in \partial_F^{-1}$. Taking a $\mathbb Z$-basis $\{1,
  \frac{1+\sqrt D}2 \}$ of $\O_F$, we see $E'\otimes_{\mathbb Z}
  \O_F = (E'\otimes 1) \oplus  (E'\otimes \frac{1+\sqrt D}2)$, and that
  $g$ is uniquely determined by (for any $e \in E$ and $x \in \O_F$)  $$
  g(e \otimes x) = \alpha(e) \otimes x + \beta(e) \otimes
  \frac{1+\sqrt D}2 x
  $$
  for some $\alpha(e), \beta(e) \in E'$, which is determined by $g$.
This implies that $\alpha$ and $\beta$ are homomorphisms from $E$
to $E'$. Let $\alpha^\vee$, $\beta^\vee$, and $g^\vee$ be dual
maps of $\alpha$, $\beta$, and $g$, then (for any $e' \in E'$ and
$y \in
\partial_F^{-1} = (\O_F)^\vee=\hbox{Hom}_{\mathbb Z} (\O_F, \mathbb Z)$)
$$
g^\vee (e'\otimes y) = \alpha^\vee(e') \otimes y + \beta^\vee
\frac{1+\sqrt D}2 y.
$$
Here we used the simple fact that with respect to the bilinear form
on $F$, $(x, y) =\tr xy$, the dual of an ideal $\mathfrak a$ is
$\mathfrak a^{-1} \partial_F^{-1}$, and the left multiplication
$l(r)$ is self-dual: $l(r)^\vee =l(r)$.

Taking $r=1$, and $x=1$, we have then for any $e \in E$
\begin{align*}
e \otimes 1 &=g^\vee \lambda_F(1) g(e\otimes 1)
 \\
  &= g^\vee (\alpha(e) \otimes 1 +  \beta(e) \otimes \frac{1+\sqrt
  D}2)
  \\
  &= \alpha^\vee \alpha(e) \otimes 1 + \beta^\vee \alpha(e) \otimes
  \frac{1+\sqrt D}2 + \alpha^\vee \beta(e) \otimes \frac{1+\sqrt D}2
  + \beta^\vee \beta(e) \otimes (\frac{1+\sqrt D}2)^2
  \\
   &= (\deg \alpha + \deg \beta \frac{D-1}4)e \otimes 1 +
   (\beta^\vee \alpha(e) + \alpha^\vee \beta(e) + \deg \beta)\otimes
    \frac{1+\sqrt D}2.
    \end{align*}
  This implies
    $$
    1 = \deg \alpha + \deg \beta \frac{D -1}4.
    $$
    So $\deg \alpha =1$ and $\deg \beta =0$. This means that
    $\alpha$ is an isomorphism, $\beta =0$, and $g =\phi(\alpha)$.
  \end{proof}

Let $\omega$  be the Hodge bundle on $\tilde{\CalM}$. Then the
rational sections of  $\omega^k$ can be identified with  meromorphic
Hilbert modular forms for $\SL_2(\O_F)$ of weight $k$. We give it
the following Petersson metric
\begin{align}\label{def:pet}
\|F(z_1,z_2)\|_{\Pet}= |F(z_1,z_2)|\big(16\pi ^2 y_1 y_2\big)^{k/2}
\end{align}
for a Hilbert modular form $F(z)$ of weight $k$. This gives a
metrized Hodge bundle $\hat{\omega}=(\omega, \|\,\|_{\Pet})$.
Strictly speaking, the metric has pre-log singularity along the
boundary $\tilde{\CalM}-\CalM$, \cite{BBK}. Since our CM cycles
never intersect with the boundary, the Faltings' height
$$
h_{\hat{\omega}}(\CM(K)) =\widehat{\deg} (\phi^* \hat{\omega})
$$
is still well-defined  where $\phi: \CM(K)  \rightarrow
\tilde{\CalM}$ is the natural map. Indeed $\phi^*\hat{\omega}$ is an
honest metrized line bundle on $\CM(K)$ as defined here. Faltings's
height for these generalized line bundle is defined in \cite{BBK}
(for schemes) which is compatible with our definition when applied
to stacks. It is proved in \cite{Ya3} that
\begin{equation} \label{neweq2.7}
h_{\hat{\omega}} (\CM(K)) = \frac{2\# \Cm(K)}{W_K} h_{\Fal}(A)
\end{equation}
for any CM abelian surface $(A, \iota, \lambda) \in \CM(\mathbb C)$.
This will be used in Section \ref{newsect7} to prove Theorem
\ref{theo1.4}.

\section{The degenerate case} \label{newsect3}

In this section, we briefly sketch a proof  of Theorem
\ref{newtheo1.3} in the degenerate case $D=1$ ($F=\mathbb Q \oplus
\mathbb Q$) which is a reformulation of Gross and Zagier's work on
singular moduli to illustrate the idea behind the proof of of
Theorem \ref{newtheo1.3}.  It also gives a new proof of the
Gross-Zagier formula on factorization of singular moduli
\cite[Theorem 1.3]{GZ2}.  Let $\mathcal M_1$ be
 the moduli stack over $\mathbb Z$ of elliptic curves,
Let $\mathcal M =\mathcal M_1 \times \mathcal M_1$ be the modular
stack over $\mathbb Z$ of pairs of elliptic curves. In this case,
$\mathcal T_1$ is the diagonal embedding of $\mathcal M_1$ into
$\mathcal M$. Let $K_i=\mathbb Q(\sqrt{d_i})$, $i =1, 2$, be
imaginary quadratic fields with fundamental discriminants $d_i <0$
and ring of integers $\O_{i}=\mathbb Z[\frac{d_i +\sqrt{d_i}}2]$,
and let $K=K_1 \oplus K_2$. For simplicity, we assume $d_i \equiv
1 \mod 4$ are prime to each other. Let $\CM(K_i)$ be the moduli
stack over $\mathbb Z$ of CM elliptic curves $(E, \iota_i)$ where
\begin{equation} \label{neweq3.1}
\iota_i: \O_{i} \subset \O_E=\End(E)
\end{equation}
such that the main involution in $\O_E$ reduces to the complex
multiplication on $\O_i$. Then $\CM(K)=\CM(K_1) \times \CM(K_2)$
is the `CM cycle' on $\mathcal M$ associated to $K$. It is easy to
see that
\begin{align}
\mathcal T_1.\CM(K)
 &=  \CM(K_1).\CM(K_2) \quad \hbox{ in } \quad \mathcal M_1 \notag
 \\
  &= \sum_{\disc[\tau_i]=d_i} \frac{4}{w_1 w_2} \log |j(\tau_1)
  -j(\tau_2)|
  \end{align}
  where $w_i=\# \O_i^*$ and $\tau_i$ are Heegner points in
  $\mathcal M_1(\mathbb C)$ of discriminant $d_i$.
So \cite[Theorem 1.3]{GZ2} can be rephrased as

\begin{theo} (Gross-Zagier)  Let the notation be as above, and  let $\tilde D =d_1 d_2$.
Then for a prime $p$, one has
 \begin{equation} \label{neweq3.2}
(\mathcal T_1.\CM(K))_p = \frac{1}{2} \sum_{ \frac{\tilde D -n^2}{4}
\in p
 \mathbb Z_{>0}}(\ord_p \frac{\tilde D-n^2}{4} +1) \beta(p,  n),
\end{equation}
where
$$ \beta(p,  n) = \prod_{l | \frac{\tilde D
-n^2}{4}} \beta_l(p,n)
$$
is given by
$$
\beta_l(p,n) =\begin{cases}
   \frac{1- \epsilon(p)^{t_p}}2 &\ff l=p,
   \\
   \frac{ 1+ (-1)^{t_l}}2 &\ff l \ne p, \epsilon(l) =-1,
  \\
   t_l + 1 &\ff   l \ne p, \epsilon(l) =1.
\end{cases}
$$
where $t_l=\ord_{l} \frac{\tilde D -n^2}{4}$, and
$$
\epsilon(l) = \begin{cases}
   (\frac{d_1}l) &\ff l \nmid d_1,
   \\
    (\frac{d_2}l) &\ff l \nmid d_2
    \end{cases}
    $$
    is as in \cite{GZ2}.
\end{theo}
\begin{proof}(sketch) The proof is a simple application of the Gross-Keating formula
\cite{GK}. A geometric point of ${\mathcal T}_1 \cap
{\CM}(K)=\mathcal T_1 \times_{\mathcal M} \CM(K)$ in
$F=\bar{\mathbb F}_p$ or $\mathbb C$ is given by a triple $(E,
\iota_1, \iota_2)$, with CM action given by (\ref{neweq3.1}).
Since $(d_1, d_2) =1$, such a point exists only when
$F=\bar{\mathbb F}_p$ with $p$ nonsplit in $K_i$ and $E$ is
supersingular. Assuming this,  $\O_E$ is a maximal order of the
unique quaternion algebra $\mathbb B$ ramified exactly at $p$ and
$\infty$. Notice that the reduced norm on $\mathbb B$ gives a
positive quadratic form on $\mathbb B$, and let $(\, , \,)$ be the
associated bilinear form.
 Let $\phi_0 =1$, $\phi_i=
\iota_i(\frac{d_i+\sqrt{d_i}}2)$, then $\iota_i$ is determined by
$\phi_i$. Let
$$
T(\phi_0, \phi_1, \phi_2) =
 \frac{1}{2} ((\phi_i, \phi_j))
 $$
 be the matrix associated to three endomorphisms $\phi_i$. Then a
 simple computation gives
 \begin{equation} \label{neweq3.4}
 T(\phi_0, \phi_1, \phi_2) =
 \begin{pmatrix} 1 & 0 &0
 \\
  \frac{d_1}2 &\frac{1}2 &0
  \\
  \frac{d_2}2 &0 &\frac{1}2
  \end{pmatrix} \diag(1, T(n))
  \begin{pmatrix} 1 &\frac{d_1}2 &\frac{d_2}2
  \\
   0 &\frac{1}2 &0
   \\
   0 &0 &\frac{1}2
   \end{pmatrix}
   \end{equation}
   with $n = 2 (\phi_1, \phi_2) -\tilde D$ and
  \begin{equation}
   T(n) = \kzxz {-d_1} {n} {n} {-d_2}.
   \end{equation}
   It is easy to see that $\frac{\tilde D-n^2}4 \in \mathbb Z_{>0}$
   (since the quadratic form is positive definite).  In general,  for
   an integer $n$ with $\frac{\tilde D-n^2}4 \in \mathbb Z_{>0}$,
   let $\tilde T(n)$ be the $3 \times 3$ matrix defined by the right
   hand side of (\ref{neweq3.4}).  If
    $\phi_i \in \O_E$ with $\phi_0=1$ satisfies $T(\phi_0, \phi_1,
    \phi_2) =\tilde T(n)$, then
    $\iota_i(\frac{d_i +\sqrt{d_i}}2) =\phi_i$
   gives actions of $\O_i$ on $E$ and thus a geometric point $(E,
   \iota_1, \iota_2)$ in the intersection. By \cite[Proposition
   5.4]{GK}, the local intersection  index $i_p(E, \iota_1,
   \iota_2)$ of $\mathcal T_1$ and $\CM(K)$ at $(E, \iota_1, \iota_2)$ depends only on $T(\phi_0, \phi_1, \phi_2)$
   and is given by (see Theorem \ref{propI0.4} and its proof for
   detail)
   \begin{equation}
i_p(E, \iota_1,
   \iota_2)= \frac{1}2 (\ord_p \frac{\tilde D-n^2}4 +1).
   \end{equation}
So
$$
(\mathcal T_1.\CM(K))_p =\frac{\log p}2 \sum_{\frac{\tilde D
-n^2}{4} \in \mathbb Z_{>0} } (\ord_p \frac{\tilde D-n^2}4 +1)
\sum_{E \, s.s.} \frac{R'(\O_E, \tilde T(n))}{\# \O_E^*}
$$
with
$$
R'(\O_E, \tilde T(n)) =\#\{ \phi_1, \phi_2 \in \O_E:\, T(1, \phi_1,
\phi_2) =\tilde T(n)\}.
$$
The summation is over isomorphic  classes  of all supersingular
elliptic curves over $\bar{\mathbb F}_p$. Next, notice that for
two supersingular elliptic curves $E_1$ and $E_2$, $\Hom(E_1,
E_2)$ is a quadratic lattice in $\mathbb B$, and they  are in the
same genus (as $E_i$ changes). Simple argument together with
 \cite[Corollary 6.23, Proposition 6.25]{GK} (see Section \ref{sect3} for
 detail) gives
\begin{align*}
\sum_{E \, s.s.} \frac{R'(\O_E, \tilde T(n))}{\# \O_E^*} &=\sum_{E
\, s.s.} \frac{R(\O_E, \tilde T(n))}{\# \O_E^* \# \O_E^*}
 \\
 &=\sum_{E_1, E_2\, s.s.} \frac{R(Hom(E_1, E_2), \tilde T(n))}{\#
 \O_{E_1}^* \#\O_{E_2}^*}
 \\
  &=\beta(p, n).
  \end{align*}
Here
$$
R(L, \tilde T(n)) =\# \{ \phi_1, \phi_2, \phi_3 \in L:\, T(\phi_1,
\phi_2, \phi_3) =\tilde T(n)\}
$$
is the representation number of representing $\tilde T(n)$ by the
quadratic lattice $L$. So $$(\mathcal T_1.\CM(K))_p = \frac{1}{2}
\sum_{ \frac{\tilde D -n^2}{4} \in
 \mathbb Z_{>0}}(\ord_p \frac{\tilde D-n^2}{4} +1) \beta(p,  n).
 $$
 Notice that $\beta_p(p, n) =0$ when $p \nmid \frac{\tilde D
 -n^2}{4}$ by the formula for $\beta_p(p, n)$. So the summation is
 really over $\frac{\tilde D -n^2}{4} \in p\mathbb Z_{>0}$.
 This proves the theorem.
\end{proof}
 In the degenerate case, it is reasonable to view $\tilde K=\mathbb
 Q(\sqrt{d_1}, \sqrt{d_2})$ as the reflex field of $K=\mathbb
 Q(\sqrt{d_1}) \oplus \mathbb Q(\sqrt{d_2})$ with respect to the `CM type' $\Phi=\{1, \sigma\}$:
 $$
  \sigma(\sqrt{d_1}, \sqrt{d_2}) = (\sqrt{d_2}, \sqrt{d_1}),\quad
 \sigma(\sqrt{d_2}, \sqrt{d_1})= (-\sqrt{d_1}, \sqrt{d_2}).
 $$
 $\tilde K$ has  real quadratic
 subfield $\tilde F=\mathbb Q(\sqrt{\tilde
 D})$ with $\tilde D =d_2 d_2$.
Using this convention, one can define $b_m(p)$ and $b_m$ as in
Conjecture \ref{conj}. We leave it to the reader to check that
$$
b_1(p) =\sum_{ \frac{\tilde D -n^2}{4} \in p
 \mathbb Z_{>0}}(\ord_p \frac{\tilde D-n^2}{4} +1) \beta(p,  n),
$$
and thus $ \mathcal T_1.\CM(K) = \frac{1}2 b_1$. This verifies
Conjecture \ref{conj} for the degenerate case $D=1$.

\section{Local Intersection  Indices} \label{sect2}

\begin{lem} \label{lemI2.1}  Let $F=\mathbb Q(\sqrt D)$ be a real quadratic field
with $D \equiv 1 \mod 4$ prime. Let $\Delta \in \O_F$ be totally
negative and let  $\tilde D =\Delta \Delta'$. Let $n$ be an
integer $0 < n < \sqrt{\tilde D}$ with $\frac{\tilde D -n^2}{D}
\in \mathbb Z_{>0}$.

(a) \quad When $D \nmid n$, there is a unique sign  $\mu =\mu(n)=
\pm 1$ and a unique positive definite integral matrix $T(\mu n)
=\kzxz {a} {b} {b} {c} \in \hbox{Sym}_2(\mathbb Z)$ such that
\begin{align}
\det T(\mu n) &= \frac{\tilde D -n^2}{D}, \label{eqI2.1}
\\
\Delta &= \frac{2 \mu n -Dc - (2b+Dc) \sqrt D}2. \label{eqI2.2}
\end{align}
Moreover, one has
\begin{equation}
a + Db + \frac{D^2-D}4 c =-\mu n. \label{eqI2.3}
\end{equation}

(b) \quad When $D |n$, for each $\mu =\pm 1$, there is a unique
positive definite integral matrix $T(\mu n) =\kzxz {a} {b} {b} {c}
\in \hbox{Sym}_2(\mathbb Z)$ such that (\ref{eqI2.1}) and
(\ref{eqI2.2}) hold. In each case, $(\ref{eqI2.3})$ holds.
\end{lem}
\begin{proof} Write $\Delta =\frac{u+v\sqrt D}2$, then
$u^2 -v^2 D =4 \tilde D$, and so
$$
u^2 \equiv 4 \tilde D \mod D \equiv 4 n^2 \mod D.
$$
This implies $D| (u-2n)(u+2n)$.

(a) \quad  Since $D \nmid 4n$ is prime, there is thus a unique $\mu
=\pm 1$ and unique integer $c$ such that
$$
u =2 \mu n -D c.
$$
Since $\Delta$ is totally negative, $u <0$. So $ u^2 \ge 4 \tilde D>
4n^2 $, and so $u < 2 \mu n$, and $c >0$. $(\ref{eqI2.2})$ also
gives $b =\frac{-v -Dc}2=\frac{u-v}2 + \mu n \in \mathbb Z$. Next,
(\ref{eqI2.1}) gives a unique $a \in \mathbb Q_{>0}$, and $T(\mu
n)>0$. We now verify that $a$ is an integer by showing that it
satisfies (\ref{eqI2.3}). The equation (\ref{eqI2.1}) gives
$$
4D ac - 4D b^2 = 4\tilde D -4 n^2 = -v^2 D -2 u D c -D^2c^2.
$$
So
\begin{align*}
4 a c &= -D c( 4b+Dc) -2(2 \mu n -Dc) c -D c^2
\\
 &=-4 D b c -D^2 c^2 + D c^2 -4 \mu n,
 \end{align*}
 and so
 $$
 a + Db +\frac{D^2 -D}4 c =- \mu n
 $$
 as claimed in (\ref{eqI2.3}).

(b) \quad  When $D|n$, $D |u$. So for each $\mu =\pm 1$, there is
a unique integer $n$ such that $u = 2 \mu n -D c$. Everything else
is the same as in (a).
\end{proof}

\begin{rem} \label{remI2.2}  Throughout this paper, the sum $\sum_{\mu}$ means either
$\sum_{ \mu =\pm 1}$ when $D|n$ or the unique term $\mu$ satisfying
the condition in Lemma \ref{lemI2.1} when $D \nmid n$.
\end{rem}

Let $E$ be a supersingular elliptic curve over $k=\bar{\mathbb
F}_p$. Then $\O_E =\End(E)$ is a maximal order of the unique
quaternion algebra $\mathbb B$ ramified exactly at $p$ and $\infty$.
Let
\begin{equation}
L_E= \{ x \in \mathbb Z + 2 \O_E:\, \tr x =0 \}
\end{equation}
be the so-called Gross lattice with quadratic form $Q(x) = x \bar
x=-x^2$, where $x \mapsto \bar x$ is the main involution of
$\mathbb B$. The reduced norm gives a quadratic form on $\mathbb
B$. For $x_1, x_2, \cdots, x_n \in \mathbb B$, we define
\begin{equation}
T(x_1, x_2, \cdots, x_n) =\frac{1}2 \left( (x_i, x_j) \right) \in
\hbox{Sym}_n(\mathbb Q).
\end{equation}

\begin{prop}  \label{propI2.2} Let the notation and assumption be as in Theorem
\ref{maintheo}. Let $p$ be a prime  and $E$ be a supersingular
elliptic curve over $\bar{\mathbb F}_p$ with endomorphism ring
$\O_E$. Then there is a one-to-one correspondence among the
following three sets.

(1) \quad The set $I(E)$ of ring embeddings $\iota: \O_K
\hookrightarrow \End_{\O_F}(E\otimes \O_F) = \O_E\otimes \O_F$
satisfying

(a) \quad $\iota(a) =1 \otimes a$ for $a \in \O_F$, and

(b) \quad the main involution in $\O_E$ induces the complex
conjugation on $\O_K$ via $\iota$.

(2) \quad The set $\mathbb T(E)$  of $(\delta, \beta) \in L_E^2$
such that $T(\delta, \beta) =T(\mu n)$ for some integer $0 < n
<\sqrt{\tilde D}$ such that $\frac{\tilde D -n^2}{4D} \in p\mathbb
Z_{>0}$ and a unique $\mu = \pm 1$.

(3) \quad The set $\tilde{\mathbb T}(E)$ of $(\alpha_0, \beta_0) \in
\O_E^2$ such that $T(1, \alpha_0, \beta_0) =\tilde T(\mu n)$ for
some integer $0 < n <\sqrt{\tilde D}$ such that $\frac{\tilde D
-n^2}{4D} \in p\mathbb Z_{>0}$ and a unique $\mu = \pm 1$. Here
$$
\tilde T =\begin{pmatrix} 1 &0 &0 \\ \frac{w_0}2 &\frac{1}2 &0 \\
\frac{w_1}2 &0 &\frac{1}2\end{pmatrix} \diag(1, T)
\begin{pmatrix} 1 &\frac{w_1}2 &\frac{w_1}2 \\ 0 &\frac{1}2 &0 \\
0 &0 &\frac{1}2\end{pmatrix} =\begin{pmatrix}
 1&\frac{w_1}2 &\frac{w_1}2 \\
 \frac{w_0}2 & \frac{1}4 (a+w_0^2) &\frac{1}4 (b+ w_0 w_1)
 \\
 \frac{w_1}2  &\frac{1}4 (b+ w_0 w_1) &\frac{1}4 (c +w_1^2)
 \end{pmatrix}
 $$
for $T=T(\mu n)$. Here $w =w_0 + w_1 \frac{D+\sqrt D}2$ is given
in (\ref{eqOK}).

 The correspondences are determined by
 \begin{align} \label{eqcorr1}
 \iota(\frac{w+\sqrt\Delta}2 )&= \alpha_0 + \beta_0 \frac{D +\sqrt
 D}2,
 \\
  \iota(\sqrt\Delta) &=\delta + \beta \frac{D+\sqrt D}2,
  \label{eqcorr2}
  \\
  \delta &=2 \alpha_0 - w_0, \quad  \beta =2 \beta_0 -w_1.
  \label{eqcorr3}
  \end{align}

\end{prop}
\begin{proof} Given an embedding $\iota \in I(E)$, we define $\alpha_0$, $\beta_0$,
$\delta$ and $\beta$ by (\ref{eqcorr1}) and  (\ref{eqcorr2}). They
satisfy (\ref{eqcorr3}), and $(\delta, \beta) \in L_E^2$. Write
$T(\delta, \beta) =\kzxz {a} {b} {b} {c} $ with $a =\frac{1}2
(\delta, \delta) =-\delta^2$, $b =\frac{1}2(\delta, \beta)$, and $c
=\frac{1}2(\beta, \beta)= -\beta^2$. First,
\begin{align*}
\Delta &=\iota(\Delta) =\iota(\sqrt\Delta)^2 =(\delta+\frac{D}2
\beta)^2 -(\delta +\frac{D}2 \beta, \frac{1}2 \beta ) \sqrt D
\\
 &= -a -Db -\frac{D^2+D}4 c - (b+\frac{1}2 Dc) \sqrt D.
 \end{align*}
 We define $n>0$ and $\mu =\pm 1$ by
 $$
 -\mu n = a + Db +\frac{D^2 -D}4 c.
 $$
 Then
 $$
 \Delta = \frac{2 \mu n-Dc -(2b +Dc)\sqrt D}2
 $$
satisfying (\ref{eqI2.2}) in Lemma \ref{lemI2.1}. Now a simple
calculation using $\tilde D =\Delta \Delta'$ gives
$$
\det T(\delta, \beta) =ac -b^2 =\frac{\tilde D -n^2}D
$$
satisfies (\ref{eqI2.1}). So $T(\delta, \beta) =T(\mu n)$ for a
unique $n$ satisfying the conditions in Lemma \ref{lemI2.1}. To show
$p |\det T(\mu n)$, let
$$
\gamma = (\delta, \beta) + 2 \delta \beta \in L_E.
$$ Then
$$(\delta, \gamma) =(\beta, \gamma)
=0, \quad (\gamma, \gamma) =2 (\delta, \delta) (\beta, \beta)-2
(\delta, \beta)^2=8 \det T(\mu n) . $$ So the determinant of $\{
\delta, \beta,\gamma \}$ is
$$
\det T(\delta, \beta, \gamma)=\det \diag( T(\mu n),  4 \det T(\mu
n))
 = 4 \det T(\mu n)^2.
    $$
Since $L_E$ has  determinant  $4p^2$, we have thus $p|\det T(\mu
n)$. To show $4|\det T(\mu n)$, it is easier to look at $\tilde
T(\mu n) \in \Sym_3(\mathbb Z)^{\vee}$ (since $\alpha_0, \beta_0 \in
\O_E$). It implies that
\begin{equation} \label{eqI2.9}
a \equiv -w_0^2 \mod 4, \quad b \equiv -w_0 w_1 \mod 2, \quad c
\equiv -w_1^2 \mod 4.
\end{equation}
So $\det T(\mu n) = ac -b^2 \equiv  0\mod 4$, and therefore
$(\delta, \beta)\in \mathbb T(E)$. A simple linear algebra
calculation shows that $(\alpha_0, \beta_0) \in \tilde{\mathbb
T}(E)$.

Next, we assume that $(\delta , \beta) \in \mathbb T(E)$. Define
$\iota$ and $(\alpha_0, \beta_0)$ by (\ref{eqcorr2}) and
(\ref{eqcorr3}). The above calculation gives
$$
(\delta +\beta \frac{D+\sqrt D}2)^2 =\Delta,
$$
so $\iota$ gives an embedding from $K$ into $\mathbb B \otimes \O_F$
satisfying the conditions in (1) once we verify $\iota(\O_K) \subset
\O_E \otimes \O_F$, which is equivalent to $\alpha_0, \beta_0 \in
\O_E$. Write
$$
\delta =-u_0 + 2 \alpha_1, \quad \beta=-u_1 + 2 \beta_1, \quad u=u_0
+ u_1 \frac{D+\sqrt D}2
$$
with $u_i \in \mathbb Z$, $\alpha_1, \beta_1\in \O_E$. Then
$$
\iota(\frac{u +\sqrt\Delta}2) = \alpha_1 + \beta_1\frac{D+\sqrt D}2
\in \O_E \otimes \O_F.
$$
This implies that $\frac{u +\sqrt\Delta}2 \in \O_K$. On the other
hand, $\frac{w+\sqrt\Delta}2 \in \O_K$.  So  $\frac{u-w}2
 \in  \O_F$, i.e., $\frac{w_i -u_i}2 \in \mathbb Z$, and
 $$
 \alpha_0 =\alpha_1 +\frac{w_0 -u_1}2 \in \O_E, \quad
 \beta_0 =\beta_1 +\frac{w_1 -u_1}2 \in \O_E
 $$
 as claimed. So $(\alpha_0, \beta_0) \in \tilde{\mathbb T}(E)$ and
 $\iota \in I(E)$. Finally, if $(\alpha_0,
 \beta_0) \in \tilde{\mathbb T}(E)$, it is easy to check that  $(\delta, \beta) \in \mathbb
 T(E)$.
\end{proof}

The proof also gives the following interesting fact. In particular,
Corollary \ref{cor1.3} is true.

\begin{cor} Let $K$ be a  non-biquadratic quartic CM number field with real
quadratic subfield $F =\mathbb Q(\sqrt D)$ where $D$ does not need
to be a prime. If $\O_K$ is a free $\O_F$-module as in
(\ref{eqOK}) with $\tilde D =\Delta \Delta' < 8D$ (not necessarily
square free or odd). There is no elliptic curve $E$ such that
$E\otimes \O_F$ has an $\O_K$-action whose restriction to $\O_F$
coincides with the natural action of $\O_F$ on $E\otimes \O_F$.
\end{cor}
\begin{proof} If such an CM action exists, $E$ has  to be a supersingular elliptic curve over $\bar{\mathbb F}_p$ for some prime $p$.
 Let $\iota$ be the resulting embedding
$\iota: \O_K \hookrightarrow \O_E \otimes \O_F$. The main
involution of $\mathbb B$ induces an automorphism of $K$ which is
the identity on $F$. Extending this through one real embedding
$\sigma$ of $F$, we get an embedding $\iota: \mathbb C
\hookrightarrow \mathbb B \otimes \mathbb R$, which  is the
division quaternion algebra over $\mathbb R$. Then main involution
has to induces the complex conjugation of $\mathbb C$ and thus
$K$. Now the same argument as above implies that there is an
integer $n > 0$ such that $\frac{\tilde D -n^2}{4D} \in p \mathbb
Z_{\ge 0}$.  Since $\tilde D$ is not a square ($K$ is not
biquadratic), one has $\frac{\tilde D-n^2}{4D} \ge p \ge 2$, i.e.,
$\tilde D \ge 8D$, a contradiction.
\end{proof}

We are now ready to deal with local intersection indices of
$\mathcal T_1$ and $\CM(K)$ at a geometric intersection point. In
view of  Lemma \ref{newlem2.3}, we consider the fiber product

\begin{equation}
\xymatrix{
  \CM(K) \times_{\CalM} \mathcal E \ar[r]^-{\tilde{f}} \ar[d]^-{\tilde{\phi}}
  &\mathcal E
  \ar[d]^-{\phi} \cr
   \CM(K) \ar[r]^-{f} &\mathcal M \cr }.
   \end{equation}
An element in $\CM(K) \times_{\CalM} \mathcal E(S)$ is a tube $(E,
A, \iota, \lambda)$ such that $(A, \iota, \lambda)\in \CM(K)(S)$
$E\in \mathcal E(S)$ satisfying
$$
A=E\otimes \O_F, \quad  \iota|_{\O_F} =\iota_F,  \quad \lambda
=\lambda_F
$$
where $\iota_F$ and $\lambda_F$ are given in Lemma
\ref{newlem2.3}.  This is determined by
$$
\iota: \O_K  \hookrightarrow \O_E \otimes \O_F
$$
with $\iota \in I(E)$.
 So an intersection point
$x=(E\otimes \O_F, \iota_F, \lambda_F) \in  \CM(K)\cap \mathcal
T_1(S)$ is given by a pair $(E, \iota)$ with $\iota \in I(E)$.
When $S=\operatorname{Spec}(F)$ for an algebraically closed field
$F=\mathbb C$ or $\bar{\mathbb F}_p$, such an pair exists only
when $F=\bar{\mathbb F}_p$ and $E$ is supersingular. Assuming
this, and write $\CalZ= \CM(K)\cap \mathcal T_1$.
 Let $W$ be the Witt ring of
$\bar{\mathbb F}_p$, and let $\mathbb E$ be the universal lifting of
$E$ to $W[[t]]$, and let $I$ be the minimal ideal of $W[[t]]$ such
that $\iota$ can be lifted to an embedding
$$
\iota_I: \O_K \hookrightarrow  \End (\mathbb E \mod I) \otimes \O_F.
$$
Then the deformation theory implies the strictly local henselian
ring $\tilde{\O}_{\CalZ, x}$ is equal to
$$
\tilde{\O}_{\CalZ, x} = W[[t]]/I.
$$
So
\begin{equation}
i_p(\CM(K), \mathcal T_1, x) =\hbox{Length}_{W} W[[t]]/I,
\end{equation}
which we also denote by $i_p(E, \iota)$.  Therefore
 \begin{equation}
 (\CM(K).\mathcal T_1)_p =  \sum_{E s.s., \,  \iota\in I(E)}  \frac{1}{\# \O_E^*} i_p(E, \iota) \log
 p.
  \end{equation}
 Here  `s.s.' stands for supersingular elliptic curves.  Notice that $\hbox{Aut}(x) =\hbox{Aut}(E) =\O_E^*$ by Lemma
 \ref{newlem2.3}.
The local intersection index $i_p(E, \iota)$ can be computed by a
  beautiful formula of Gross and Keating \cite{GK} as follows.

  \begin{theo}  \label{propI0.4} Let the notation be as above, and let $(\delta , \beta) \in \mathbb T(E)$ be the image of $\iota \in I(E)$, and
  let $T(\mu n) =T(\delta, \beta)$ as in Proposition \ref{propI2.2}. Then
  $$
  i_p(E, \iota)
   = \frac{1}2(\ord_p \frac{\tilde D-n^2}{4D} +1)
     $$
     depends only on $n$.
     \end{theo}
\begin{proof}  Let $(\alpha_0, \beta_0) \in \tilde{\mathbb T}(n)$ be the image of $\iota$.
First notice that $I$ is also the smallest ideal of $W[[t]]$ such
that $\alpha_0$ and $\beta_0$ can be lifted to endomorphisms of
$\mathbb E \mod I$. Applying the Gross-Keating formula
\cite[Proposition 5.4]{GK} to $f_1=1$, $f_2=\alpha_0$, and
$f_3=\beta_0$, we see that $i_p(K, \iota)$ depends on the
$\GL_3(\mathbb Z_p)$-equivalence class of $\tilde T(\mu n)$ and is
given as follows.

  Let $a_0
\le a_1 \le a_2$ be the Gross-Keating invariants of the quadratic
form
\begin{equation}
Q(x + y \alpha_0 + z \beta_0)= (x, y, z)  \tilde T(\mu n)  (x, y,
z)^t
\end{equation}
defined in \cite[Section 4]{GK}. Then $i_p(E, \iota)$ equals
\begin{align*}
&\sum_{i=0}^{a_0-1} (i+1) (a_0+a_1+a_2-3i) p^i
+\sum_{i=a_0}^{(a_0+a_1-2)/2} (a_0+1)(2a_0+a_1+a_2-4i) p^i
\\
 &\qquad
  +\frac{a_0+1}2(a_2-a_1+1)p^{\frac{a_0+a_1}2}
\end{align*}
if $a_1 - a_0$ is even, and
$$
\sum_{i=0}^{a_0-1} (i+1)(a_0+a_1+a_2-3i) p^i
+\sum_{i=a_0}^{(a_0+a_1-1)/2} (a_0+1) (2a_0+a_1+a_2-4i) p^i
$$
if $a_1-a_0$ is odd.

First assume that  $p$ is odd. In this case, $\tilde T(\mu n)$ is
$\GL_3(\mathbb Z_p)$-equivalent to $\diag(1, T(\mu n))$. Notice that
$p\nmid T(\mu n)$,    $T(\mu n)$ is $\GL_2(\mathbb Z_p)$-equivalent
to $\diag(\alpha_p, \alpha_p^{-1} \det T(\mu n))$ for some $\alpha_p
\in \mathbb Z_p^*$, so $\tilde T(\mu n)$ is equivalent to $\diag(1,
\alpha_p,\alpha_p^{-1} \det T(\mu n))$.  So the Gross-Keating
invariants are $(0, 0, \ord_p \det T(\mu n))$ in this case. The
Gross-Keating formula gives
$$
i_p(E, \iota)=\frac{1}2 (\ord_p \det T(\mu n)+1) =\frac{1}2 (\ord_p
\frac{\tilde D-n^2}{4D} +1).
$$

Now we assume $p=2$. Since the quadratic form $Q$ associated to
$\tilde T(\mu n)$ is anisotropic over $\mathbb Q_2$, $\tilde T(\mu
n)$ is $\GL_3(\mathbb Z_2)$-equivalent to either
$$
\diag(\epsilon_0 2^{t_0}, 2^s \kzxz {1} {1/2} {1/2} {1}) \quad
\hbox{ or } \quad \diag(\epsilon_1 2^{t_1}, \epsilon_2
2^{t_2},\epsilon_3 2^{t_3})
$$
with $\epsilon_i \in \mathbb Z_2^*$ and $t_i, s\in \mathbb Z_{\ge
0}$. Since $\tilde T(\mu n)$ is not integral over $\mathbb Z_2$ (at
least one of $w_0$ or $w_1$ is odd), $\tilde T(\mu n)$ has to be
$\GL_3(\mathbb Z_2)$-equivalent to $ \diag(\epsilon_0 2^{t_0}, \kzxz
{1} {1/2} {1/2} {1})$. In this case, \cite[Proposition
B.4]{YaDensity2} asserts that the Gross-Keating invariants are $(0,
0, t_0)$. Since
$$
3 \epsilon_0 2^{t_0 -2} =\det \tilde T(\mu n) =\frac{1}{16} \det
T(\mu n) = \frac{1}4 \frac{\tilde D -n^2}{4 D},
$$
we see $t_0 =\ord_2  \frac{\tilde D -n^2}{4 D}$. Now the
Gross-Keating formula gives the desired  formula for $p=2$.
\end{proof}

We remark that when $p \ne 2$, the ideal $I$ is also the minimal
ideal such that $\delta$ and $\beta$ can be lifted to endomorphisms
of $\mathbb E \mod I$. It is not true for $p=2$.  In summary, we
have our first main formula.

\begin{theo} \label{newtheo4.7} Let the notation and assumption be as in Theorem
\ref{maintheo}. Then
$$
 \CM(K).\mathcal T_1=\frac{1}2 \sum_p \log p \sum_{0<n<\sqrt{\tilde D},\, \frac{\tilde D -n^2}{4D} \in p\mathbb Z_{>0}}
 (\ord_p \frac{\tilde D -n^2}{4D}+1)\sum_\mu  \sum_{E\,  s.s} \frac{R(L_E, T(\mu n))}{\# \O_E^*}
$$
where   $R(L, T)$ is the number of a lattice $L$ representing $T$,
i.e.,
$$
R(L, T) = \# \left\{ x=(x_1, x_2) \in L^2:\,  T(x_1, x_2)=T
\right\}.
$$
Here the meaning of $\sum_\mu$ is given in Remark \ref{remI2.2}.
\end{theo}

\section{Local Density} \label{sect3}

It is not hard to prove that the quantity $\sum_{E s.s} \frac{R(L_E,
T(\mu n))}{\# \O_E^*}$ is product of the so-called local densities.
Explicit formulae  for these local densities are given by the author
in \cite[Propositions 8.6 and 8.7]{YaDensity1}  for $p \ne 2$.  When
$p=2$, an algorithm is given in \cite{YaDensity2} and could be used
to get a formula in our case, but it is cumbersome. We give an
alternative way to compute it in  this section. We first prove

\begin{lem} \label{lemI3.1} Let $T \in \Sym_{2}(\mathbb Z)$ and let $\tilde T$ be
obtained from $T$ by the formula in Proposition \ref{propI2.2}(3).
Assume $\tilde T \in \Sym_3(\mathbb Z)^{\vee}$. Then
$$
\sum_{E \, s.s.}\frac{R(L_E, T)}{\# \O_E^*} =\sum_{E, E'\,
s.s}\frac{R(\hbox{Hom}(E, E'), \tilde  T)}{\# \O_E^* \# \O_{E'}^*}.
$$
Here $\hbox{Hom}(E, E')$ is equipped with the degree as its
quadratic form.
\end{lem}
\begin{proof} Clearly $T(\delta, \beta) =T$ if and only if
$T(1, \alpha_0, \beta_0) =\tilde T$ where $\delta, \beta, \alpha_0$,
and $\beta_0$ are related by (\ref{eqcorr3}). The condition $\tilde
T \in \Sym_3(\mathbb Z)^{\vee}$ implies that $(\delta, \beta) \in
L_E^2$ if and only if $(\alpha_0, \beta_0) \in \O_E^2$. Next if
$f_1, f_2, f_3 \in \O_E$ represents $\tilde T$, i.e., $T(f_1, f_2,
f_3) =\tilde T$, then the reduced norm $N f_1 =1$, and $f_1 \in
\O_E^*$. In such a case,  $1, \alpha_0 =f_1^{-1} f_2,
\beta_0=f_1^{-1} f_3$ represents $\tilde T$ too. So
$$
R(\O_E, \tilde T)= \# \O_E^* R(L_E, T)
$$
and
$$
\sum_{E\, s.s} \frac{R(L_E, T)}{\# \O_E^*} =\sum_{E\, s.s}
\frac{R(\O_E, \tilde T)}{(\# \O_E^*)^2}.
$$
On the other hand, if $f_1, f_2, f_3 \in \hbox{Hom}(E, E')$
represents $\tilde T$, then $\deg f_1 =1$ and $f_1$ is actually an
isomorphism. So $R(\hbox{Hom}(E, E'), \tilde T)=0$ unless $E$ and
$E'$ are isomorphic. This proves the lemma.
\end{proof}

The quantity
$$
\sum_{E, E'\, s.s}\frac{R(\hbox{Hom}(E, E'), \tilde  T)}{\# \O_E^*
\# \O_{E'}^*}
$$
is also a product of local densities and is computed by Gross and
Keating \cite[Section 6]{GK} in terms of Gross-Keating invariants
(see also \cite{We1}  and \cite{We2} for more extensive
explanation). More precisely, we have (note that $\det Q$ in
\cite{GK} is our $\det 2\tilde T$),

\begin{equation} \label{eqI3.1}  \sum_{E, E'\, s.s}\frac{R(\hbox{Hom}(E, E'), \tilde  T)}{\#
\O_E^* \# \O_{E'}^*} =\prod_{ l | 4 p\det \tilde T} \beta_l(\tilde
T).
\end{equation}
Here $\beta_p(\tilde T)$ is $0$ or $1$ depending on whether $\tilde
T$ is isotropic or anisotropic over $\mathbb Z_p$. For $l \ne p$,
$\beta_l(\tilde T)$ is zero if $\tilde T$ is anisotropic over
$\mathbb Z_l$  and is  given as follows if $\tilde T$ is isotropic
by \cite[Proposition 6.24]{GK}. Let $0 \le a_0 \le a_1 \le a_2$ be
Gross-Keating invariants of $\tilde T$, and let $\epsilon=\pm 1$ be
the Gross-Keating epsilon sign of $\tilde T$, then

If $a_0 \equiv a_1 \mod 2$ and $\epsilon=1$, we have
$$
\beta_l(\tilde T)= 2 \sum_{i=0}^{a_0-1} (i+1) l^i + 2
\sum_{i=a_0}^{(a_0 +a_1 -2)/2}(i+1) l^i + (a_0+1) (a_2-a_1+1)
l^{\frac{a_0+a_1}2}.
$$

If $a_0 \equiv a_1 \mod  2$ and $\epsilon =-1$, we have
$$
\beta_l(\tilde T)= 2 \sum_{i=0}^{a_0-1} (i+1) l^i + 2
\sum_{i=a_0}^{(a_0 +a_1 -2)/2}(i+1) l^i + (a_0+1)
l^{\frac{a_0+a_1}2}.
$$

If $a_0 \not\equiv a_1 \mod 2$, we have
$$
\beta_l(\tilde T)= 2 \sum_{i=0}^{a_0-1} (i+1) l^i + 2
\sum_{i=a_0}^{(a_0 +a_1 -1)/2}(i+1) l^i.
$$


\begin{lem}  \label{lemI3.2} Let $n$ be a positive integer such that $\frac{\tilde D
-n^2}{4D} \in \mathbb Z_{>0}$ and let $T(\mu n)$ be as in Lemma
\ref{lemI2.1}. Then

(1) \quad $T( \mu n)$ is $\GL_2(\mathbb Z)_l$-equivalent to
$\diag(\alpha_l, \alpha_l^{-1}  \det T(\mu n))$ for some $\alpha_l
\in \mathbb Z_l$.

(2) \quad   $\tilde T(\mu n)$ is isotropic if and only if
$(-\alpha_l, l)_l^{t_l} =1$ where $t_l =\ord_l \frac{\tilde D
-n^2}{4D}$.
\end{lem}
\begin{proof}  (1) follows from  the fact $l \nmid a$ or $l
\nmid c$. By (1), $\tilde T(\mu n)$ is $\GL_3(\mathbb
Q_l)$-equivalent to $\diag(1, \alpha_l, \alpha_l^{-1} \det T(\mu
n))$. So its Hasse invariant is $(\alpha_l, -\det T(\mu n))_l$. By
\cite[Chapter 4, Theorem 6]{Se}, $\tilde T(\mu n)$ is isotropic over
$\mathbb Z_l$ if and only if
$$
(\alpha_l, -\det T(\mu n))_l=(-1, -\det \tilde T(\mu n))_l,
$$
i.e.,
$$
(-\alpha_l, - \frac{\tilde D -n^2}{4D})_l =1.
$$
 When $l\ne 2$,
$$
(-\alpha_l, -\frac{\tilde D -n^2}{4D})_l =(-\alpha_l, l)_l^{t_l}.
$$
When $l =2$, Lemma \ref{lemI4.1} in the next section asserts that
either $a \equiv -1 \mod 4$ or $c \equiv -1 \mod 4$. So $\alpha_2
=a$ or $c$, and thus $-\alpha_2 \equiv 1 \mod 4$. This implies
$$
(-\alpha_2, - \frac{\tilde D -n^2}{4D})_2 =(-\alpha_2, 2)_2^{t_2}.
$$
just as in the odd case.
\end{proof}


\begin{prop}  \label{propI3.3} Assume that $ 0 < n <\sqrt{\tilde D}$ with
$\frac{\tilde D -n^2}{4D} \in p\mathbb Z_{>0}$. Then
$$
\sum_{E \, s.s} \frac{R(L_E, T(\mu n))}{\# \O_E^*} =\prod_{l|
\frac{\tilde D-n^2}{4D}} \beta_l(p, \mu n).
$$
Here $\beta_l(p, \mu n)$ is given as follows. Let $T(\mu n)$ be
$\GL_2(\mathbb Z_l)$-equivalent to $\diag(\alpha_l, \alpha_l^{-1}
\det T(\mu n))$ over $\mathbb Z_l$ with $\alpha_l \in \mathbb
Z_l^*$, and write $t_l =\ord_l \frac{\tilde D-n^2}{4D}$. Then

$$ \beta_l(p,\mu n) =\begin{cases}
 \frac{1- (-\alpha_p, p)_p^{t_p}}2 &\ff l=p,
 \\
 \frac{1 + (-1)^{t_l}}2 &\ff l \ne p, \hbox{ and  } (-\alpha_l, l)_l =-1,
 \\
 t_l + 1 &\ff l \ne p, \hbox{ and } (-\alpha_l, l)_l =1.
 \end{cases}
$$

\end{prop}
\begin{proof} By Lemma \ref{lemI3.1} and (\ref{eqI3.1}), it
suffices to verify the formula for $\beta_l(p, n)$. The case $l=p$
follows from Lemma \ref{lemI3.2}.

 When $l \nmid 2p$, $\tilde
T(\mu n)$ is $\GL_3(\mathbb Z_l)$-equivalent to $\diag(1, T)$ and
thus to $\diag(1, \alpha_l, \alpha_l^{-1} \det T(\mu n))$. When it
is isotropic, its Gross-Keating epsilon sign is $(-\alpha_l, l)_l$
by definition \cite[Section 3]{GK}. So its Gross-Keating invariants
are $(0, 0, t_l)$,  and when it is isotropic,  its Gross-Keating
epsilon sign is $(-\alpha_l, l)_l$ by definition \cite[Section
3]{GK}. Now the formula follows from Lemma \ref{lemI3.2} and the
Gross-Keating formula described before the lemma.

Now we assume  $l=2 \ne p$. Since $\tilde T(\mu n) \notin
\Sym_3(\mathbb Z_2)$, $\tilde T(\mu n)$ is $\GL_3(\mathbb
Z_2)$-equivalent to
$$
\diag(\epsilon 2^{t_2}, \kzxz {A} {1/2} {1/2} {A}),  \quad  A =0, 1.
$$
It is isotropic if and only if $A =0$ or $A=1$ and $t_2$ is even.
 In each case, the Gross-Keating invariants are $(0,
0, t_2)$ by \cite[Proposition B.4]{YaDensity2}. In the isotropic
case, the Gross-Keating epsilon sign is $1$ if $A=0$ and $-1$ if
$A=1$ by the same proposition. We claim that $A=0$ if and only if
$(-\alpha_2, 2)_2 =1$, i.e., $\alpha_2 \equiv \pm 1 \mod 8$, i.e.,
the Gross-Keating epsilon sign of $\tilde T(\mu n)$ is again
$(-\alpha_2, 2)_2$. Indeed, Lemma \ref{lemI4.1} implies that $a
\equiv 3 \mod 4$ or $c \equiv 3 \mod 4$. Assume without loss of
generality $a \equiv 3 \mod 4$. In this case $w_0 \equiv 1 \mod 2$
and  we can take $\alpha_2 =a$. It is easy to see that $\tilde T(\mu
n)$ is $\mathbb Z_2$-equivalent to
$$
T = \begin{pmatrix} 1 &\frac{1}2 & \frac{w_1}2
   \\
    \frac{1}2  &\alpha  &\frac{1}4 (b+ w_1)
    \\
     \frac{w_1}2 &\frac{1}4(b+w_1) & \frac{1}4(c+w_1^2)
     \end{pmatrix}
     $$
     with $\alpha =\frac{1}4(a+1) \in \mathbb Z_2$.

If $(-a, 2)_2 =1$, i.e., $a \equiv 7 \mod 8$, and so $\alpha
=\epsilon 2^r$ for some $r \ge 1$ and $\epsilon \in \mathbb Z_2^*$.
Let $\beta_1$ and $\beta_2$ are roots of $x^2 + x +\alpha =0$ with
$\beta_1 \in \mathbb Z_2^*$ and $\beta_2 \in 2^r \mathbb Z_2^*$, and
let $L=\oplus \mathbb Z_2 e_i$ be the lattice of $T$. Let
$$
f_1 =e_2 + \beta_1 e_1, \quad f_2 =e_2 + \beta_2 e_2.
$$
Then it is easy to check that $(f_1, f_2) =-1+4 \alpha$ and $(f_1,
f_1) =(f_2, f_2)=0$. This implies that $L=(\mathbb Z_2 f_1 + \mathbb
Z_2 f_2) \oplus \mathbb Z_2 f_3$ for some $f_3 \in L$.  So $T$ and
thus $\tilde T(\mu n)$ is $\mathbb Z_2$-equivalent to
$$\diag( \kzxz {0} {1/2} {1/2} {0}, \epsilon_1 2^{t_2}).
$$

If $(-a, 2)_2=-1$, then $a \equiv 3 \mod 8$ and $\alpha \in \mathbb
Z_2^*$. In this case, it is easy to check  by calculation that $L
=\mathbb Z_2 e_1 \oplus \mathbb Z_2 e_2 \oplus \mathbb Z_2 f_3$ for
some $f_3 \in L$ perpendicular to $e_1$ and $e_2$, and its quadratic
form is
$$
Q(x e_1 +y e_2 +z f_3) = x^2 + xy + \alpha y^2 + d z^2
$$
with $d =Q(f_3)$, and is thus $\mathbb Z_2$-equivalent to $x^2 + xy
+y^2 + d_1 z^2$. So $A=1$. This proves the claim. The claim implies
that the formula is also true for $l =2$.

\end{proof}

{\bf Proof of Theorem \ref{newtheo1.3}}: Now Theorem
\ref{newtheo1.3} is clear from Theorem \ref{newtheo4.7} and
Proposition \ref{propI3.3}

\section{Computing $b_1(p)$ and Proof of Theorem \ref{maintheo}}
\label{sect4}

The formula for $b_1(p)$ is known to be independent of the choice of
the CM type $\Phi$. We choose $\Phi =\{1, \sigma\}$ with $\sigma
(\sqrt{\Delta}) =\sqrt{\Delta'}$. Then $\tilde K =\tilde F
(\sqrt{\tilde\Delta})$ with
\begin{equation}
\tilde\Delta =(\sqrt\Delta + \sqrt{\Delta'})^2 = \Delta + \Delta' -2
\sqrt{\tilde D}.
\end{equation}
For an integer $0 < n <\sqrt{\tilde D}$, let $\mu$ and $T(\mu n)
=\kzxz {a} {b} {b} {c} $ be as in Lemma \ref{lemI2.1}. Then we have
by Lemma \ref{lemI2.1}
\begin{equation}
\tilde\Delta = 2\mu n -Dc  - 2 \sqrt{\tilde D}
\end{equation}

\begin{lem}  \label{lemI4.1} Assume the condition $(\clubsuit)$. For an integer $0 < n < \sqrt{\tilde D}$
with $\frac{\tilde D -n^2}{4D}\in \mathbb Z_{>0}$. let $\mu$ and
$T(\mu n)=\kzxz {a} {b} {b} {c} $ be as in Lemma \ref{lemI2.1}. Then

(1) \quad One has  $\frac{\mu n +\sqrt{\tilde D}}{2D} \in d_{\tilde
K/\tilde F}^{-1}-\O_{\tilde F}$ and $ \frac{-\mu n +\sqrt{\tilde
D}}{2D} \in d_{\tilde K/\tilde F}^{-1, \prime}-\O_{\tilde F}$. Here
$'$ stands for  the Galois conjugation in $\tilde F$.

(2) \quad For any prime $p|\frac{\tilde D-n^2}{D}$, $p\nmid a$ or $p
\nmid c$.

(3) \quad Exactly one of $a$ and $c$ is $0 \mod 4$ and the other is
$-1 \mod 4$.

\end{lem}
\begin{proof}  (1): When  $D|n$,  one has  $D|\tilde D$, and $d_{\tilde K/\tilde F} =d_{\tilde K/\tilde F}'$ and the claim is clear.
When $D \nmid n$, one has  $\frac{\pm n + \sqrt{\tilde D}}{2D}
\notin \O_{\tilde F}$, and
$$
\frac{\tilde D -n^2}{4} =\frac{n+\sqrt{\tilde D}}2 \cdot
\frac{-n+\sqrt{\tilde D}}2  \in D\mathbb Z=d_{\tilde K/\tilde F}
d_{\tilde K/\tilde F}'. $$
 So there is a unique $\nu =\pm
1$ such that
$$\frac{\nu n +\sqrt{\tilde D}}2 \in d_{\tilde K/\tilde
F}', \quad \frac{-\nu n +\sqrt{\tilde D}}2 \in d_{\tilde K/\tilde
F}. $$ On the other hand, $\tilde\Delta \in d_{\tilde K/\tilde F}$
implies $\mu n-\sqrt{\tilde D} \in d_{\tilde K/\tilde F}$. So $(\mu
-\nu) n \in d_{\tilde K/\tilde F}$ and thus $(\mu -\nu) n \equiv 0
\mod D$. So $\mu =\nu$. Now it is easy to see
$$
\frac{\mu n +\sqrt{\tilde D}}{2D} \in d_{\tilde K/\tilde
F}^{-1}-\O_{\tilde F}, \quad \frac{-\mu n +\sqrt{\tilde D}}{2D}
\notin d_{\tilde K/\tilde F}^{-1,\prime}-\O_{\tilde F}.
$$

(2): If $p|a, c$, then $p|ac-b^2 =\frac{\tilde D-n^2}{D}$ implies
$p|b$, and thus $$ p|n =-\nu (a +D b +\frac{D^2-D}4 c).$$ This
implies $p| \tilde D$.  But this causes a contradiction: $p^2 |ac
-b^2 =\frac{\tilde D-n^2}{D}$.

(3): Since $K/F$ is unramified at primes of $F$ over $2$ under the
condition $(\clubsuit)$, there are integers $x$ and $y$, not both
even, such that
$$
\Delta \equiv (x + y \frac{1+\sqrt D}2)^2 \mod 4.
$$
By Lemma \ref{lemI2.1}, this implies
$$
-a =(D+1)b -\frac{D^2-D}4 c -D c \frac{1+\sqrt D}2 \equiv  x^2 + y^2
\frac{D-1}4 +(y^2 + 2 xy ) \frac{1+\sqrt D}2 \mod 4.
$$
So (since $D \equiv 1 \mod 4$)
\begin{align}
a +x^2 + 2b + \frac{D-1}4 (c+y^2) &\equiv 0 \mod 4, \label{eqI4.3}
\\
2x y + y^2 +c &\equiv 0 \mod 4. \label{eqI4.4}
\end{align}

When  $x$ is even, $y$ has to be odd.  So $c \equiv -1 \mod 4$ and
$a$ is even.
 Since
\begin{equation} \label{eqI4.5}
ac-b^2 =\det T(\mu n) =\frac{\tilde D-n^2}{D} \equiv 0 \mod 4,
\end{equation}
 one has $a  \equiv 0 \mod 4$ and $b \equiv 0 \mod 2$.

When $x$ is odd and $y$ is even, (\ref{eqI4.3}), (\ref{eqI4.4}), and
(\ref{eqI4.5}) imply $a+1 \equiv c \equiv 0 \mod 4$.

When both $x$ and $y$ are odd, (\ref{eqI4.4}) implies  $c \equiv 1
\mod 4$.  So (\ref{eqI4.3}) implies that $a$ is odd and  thus $b$
is odd. Now  (\ref{eqI4.5})  implies $a \equiv 1 \mod 4$ and that
$b$ is odd. So (\ref{eqI4.3}) implies $D \equiv 1 \mod 8$. But
this implies that $ -\mu n =a +D b + \frac{D^2 -D}4 c$  is even,
which is impossible since $\tilde D$ is odd. So this case is
impossible.
\end{proof}

It is easy to see from the definition that $b_1(p)=0$ unless  there
is $n
>0$ with  $ \frac{\tilde D -n^2}{4D} \in p \mathbb Z_{>0}$.
  This implies in particular $p$ is split in
$\tilde F$ or $p|\tilde D$ is ramified in $\tilde F$. For a fixed $n
>0$  with $ \frac{\tilde D -n^2}{4D} \in p \mathbb Z_{>0}$, fix a sign $\mu =\pm 1$ such  that  $T(\mu n)$ exists as in Lemma \ref{lemI2.1}. In the split case, we choose
the splitting $p\O_{\tilde F} = \mathfrak p \mathfrak p'$ such that
\begin{equation}
t_{\mu n} =\frac{\mu n+ \sqrt{\tilde D}}{2D} \in \mathfrak p
d_{\tilde K/\tilde F}^{-1}.
\end{equation}
So
\begin{equation}
\ord_\mathfrak p t_{\mu n} = \ord_p \frac{\tilde D -n^2}{4D}, \quad
\ord_{\mathfrak p'} (t_{\mu n}) =0  \hbox{ or } -1.
\end{equation}
In the ramified case $p\O_{\tilde F} =\mathfrak p^2$, the above two
equations also hold (forgetting the one with $\mathfrak p'$). With
this notation, we have by (\ref{eqI1.3})
\begin{equation}
b_1(p) = \sum_{ 0< n <\sqrt{\tilde D}, \frac{\tilde D -n^2}{4D} \in
p \mathbb Z_{>0}} (\ord_p \frac{\tilde D -n^2}{4D} +1)\sum_{\mu}
b(p,\mu n)
\end{equation}
where
\begin{equation}
b(p, \mu n) = \begin{cases}
  0 &\ff \mathfrak p \hbox { is  split in } \tilde K,
  \\
  \rho(t_{\mu n} d_{\tilde K/\tilde F} \mathfrak p^{-1}) &\ff \mathfrak p \hbox { is not  split in } \tilde
  K.
  \end{cases}
  \end{equation}

Assume now that $\mathfrak p$ is not split in $\tilde K$. Notice
that
$$
\rho(t_{\mu n}d_{\tilde K/\tilde F} \mathfrak
p^{-1})=\prod_{\mathfrak l} \rho_{\mathfrak l}(t_{\mu n}d_{\tilde
K/\tilde F} \mathfrak p^{-1})
$$
where the product runs over all  prime ideals  $\mathfrak l$ of
$\tilde F$, and
\begin{equation} \label{eqI2.4}
\rho_{\mathfrak l}(t_{\mu n} d_{\tilde K/\tilde F}\mathfrak p^{-1})
 = \begin{cases}
   1    &\ff    \mathfrak l \hbox{ is ramified in  } \tilde K,
   \\
   \frac{1 +(-1)^{\ord_{\mathfrak l} (t_{\mu n} d_{\tilde K/\tilde F} \mathfrak
   p^{-1})}}2 &\ff  \mathfrak l \hbox{ is inert in  } \tilde K,
   \\
   1+\ord_{\mathfrak l} (t_{\mu n} d_{\tilde K/\tilde F }\mathfrak
   p^{-1})
    &\ff  \mathfrak l \hbox{ is split  in  } \tilde K.
    \end{cases}
    \end{equation}
Notice first that $\rho_{\mathfrak l}(t_{\mu n} d_{\tilde K/\tilde
F}\mathfrak p^{-1})=1$ unless $\ord_{\mathfrak l} (t_{\mu n}
d_{\tilde K/\tilde F}\mathfrak p^{-1}) \ge 1$. In  particular it is
$1$ unless $l | \frac{\tilde D -n^2}{4 D p}$ where $l$ is the prime
under $\mathfrak l$.  So we assume that $l | \frac{\tilde D -n^2}{4
D p}$. This implies that either $l|\tilde D$ is ramified in $\tilde
F$ or  $l=\mathfrak l \mathfrak l'$ is split in $\tilde F$. In the
split case,   we  again choose the splitting $l \O_{\tilde F}
=\mathfrak l \mathfrak l'$ so that
\begin{equation} \label{eqI2.5}
\ord_{\mathfrak l} (t_{\mu n} d_{\tilde K/\tilde F}\mathfrak p^{-1})
=\ord_l \frac{\tilde D -n^2}{4 D p}=\ord_{l} \frac{\det T(\mu
n)}{4p}, \quad \ord_{\mathfrak l'} (t_{\mu n} d_{\tilde K/\tilde
F}\mathfrak p^{-1})=0.
\end{equation}

\begin{lem}  \label{lemI4.2} Let the notation be as above. Assume $\mathfrak l \ne d_{\tilde K/\tilde F}$ and  that $T(\mu n)$ is
$\GL_2(\mathbb Z_l)$-equivalent to $\diag(\alpha_l, \alpha_l^{-1}
\det T(\mu n))$ with $\alpha_l \in \mathbb Z_l^*$. Then $\tilde
K/\tilde F$ is split (resp. inert) at $\mathfrak l$ if and only if
$(-\alpha_l, l)_l=1$ (resp. $-1$).
\end{lem}
\begin{proof} Since  $\mathfrak l \ne d_{\tilde K/\tilde F}$,
 $\mathfrak l$ is either split or inert in $\tilde K$.
We first take care of a  very special case $\mathfrak l=d_{\tilde
K/\tilde F}' \ne d_{\tilde K/\tilde F}$. In this case, $l=D \nmid
n$. So
$$\tilde\Delta =2 \mu n -Dc -2 \sqrt{\tilde D} \equiv 2 \mu n -2
\sqrt{\tilde D} \equiv 4\mu n  \equiv -4 a \mod \mathfrak l
$$
So $l\nmid a$, and $\tilde K=\tilde F(\sqrt{\tilde\Delta})$ is split
at $\mathfrak l$ if and only if $(-a, l)_l=1$. In this case
$\alpha_l =a$.

  Now we can assume $l \ne D$. We divide  the proof into three cases and more subcases.

  {\bf Case 1}: First we assume that $l | \tilde D$ is ramified in
  $\tilde F$. In this case, $l\O_{\tilde F} =\mathfrak l^2$,
  $l|n$, and
  $$
  \tilde\Delta = 2 \mu n -D c - 2 \sqrt{\tilde D} \equiv - D c \mod
  \mathfrak l.
  $$
  Since $\ord_l \frac{\tilde D -n^2}{4D} =1$, it is easy to check
  $l\nmid c$ and thus $\tilde\Delta \not\equiv 0 \mod \mathfrak l$.
  Since $ \tilde\Delta {\tilde\Delta}' = D v^2$, one has
  $$
  \tilde\Delta^2  \equiv D v^2 \mod \mathfrak l
  $$
  and
$$
1 \equiv - \tilde\Delta c v^2 \mod \mathfrak l.
$$
So $\tilde K =\tilde F(\sqrt{\tilde\Delta})$ is split at $\mathfrak
l$ if and only if $(-c, l)_l =1$. Notice that $\alpha_l= c$ in the
case.

{\bf Case 2}: Next we assume that $l \nmid 2 \tilde D$. So $l$ is
split as discussed above.  In this case,   either $\mathfrak l \nmid
\tilde\Delta$ or $\mathfrak l \nmid (\tilde\Delta)'$.

{\bf Subcase 1}: We first assume $\mathfrak l \nmid \tilde\Delta$.
Since $t_{\mu n} \in \mathfrak l$, one has
$$
\tilde\Delta = 2 \mu n -Dc -2 \sqrt{\tilde D} \equiv 4 \mu n -Dc
\mod 4\mathfrak  l \equiv -\alpha \mod 4\mathfrak l
$$
with $\alpha =4a +4 Db + D^2 c$. So over $\tilde F_\mathfrak l
=\mathbb Q_l$, one has
$$
(\tilde\Delta, l)_l = (-\alpha, l)_l.
$$
So  $(-\alpha, l)_l =1$ (resp. $-1$) if and only if $\tilde K/\tilde
F$ is split (resp. inert) at $\mathfrak l$. On the other hand,
$$
\kzxz {2} {D} {1} {\frac{D+1}2} T(\mu n) \kzxz {2} {1} {D}
{\frac{D+1}2} = \kzxz {\alpha } {*} {*} {*},
$$
 $T(\mu n)$ is $\GL_2(\mathbb Z_l)$-equivalent to $\diag(\alpha,
\alpha^{-1} \det T(\mu n))$.

{\bf Subcase 2}: We now assume $\mathfrak l \nmid {\tilde\Delta}'$.
Since $\tilde\Delta {\tilde\Delta}' =D v^2$ for some integer $v$, we
see that $$(\tilde\Delta, l)_l=({\tilde\Delta}', l)_l (D, l)_l.$$ On
the other hand,
$$
{\tilde\Delta}' =2\mu n-Dc +2\sqrt{\tilde D} \equiv -Dc \mod
4\mathfrak l, $$ one sees that $(\tilde\Delta, l)_l=(-c, l)_l$ and
$l \nmid c$. Since $T(\mu n)$ is $\GL_2(\mathbb Z_l)$-equivalent to
$\diag(c, c^{-1}\det T(\mu n))$ in this case,  the lemma is true
too.

{\bf Case 3}: Finally, we deal with the case  $l= 2$. When $2\nmid
c$,  the same argument as in Subcase 2 above gives the lemma. When
$2|c$, the situation is more complicated and technical. First, $ac
-b^2 =\frac{\tilde D -n^2}{D} \equiv 0 \mod 8$ implies that $4| c$
and $2|b$. We choose the splitting $2 =\mathfrak l \mathfrak l'$ as
in (\ref{eqI2.5}).

{\bf Subcase 1}: \quad  If $c=8 c_1 \equiv 0 \mod 8$, then $b=4 b_1
\equiv 4$ and
\begin{align*}
-\frac{\tilde\Delta}4 &= -\mu n +2 D c_1 +\frac{\mu n +\sqrt{\tilde
D}}2
 \\
  &= a + 4 D b_1 + 2 D^2 c_1 + \frac{\mu n +\sqrt{\tilde D}}2
  \\
  &\not\equiv 0 \mod \mathfrak l.
  \end{align*}
Notice that
\begin{equation}\label{eqI4.9}
(\frac{x+y\sqrt{\tilde D}}2)^2 = (\frac{x -\mu n y}2)^2 + xy
\frac{\mu n +\sqrt{\tilde D}}2 + y^2 \frac{\tilde D -n^2}{4}.
\end{equation}
So $\tilde K =\tilde F(\sqrt{\tilde{\Delta}})$ is split at
$\mathfrak l$ if and only if there is $\frac{x+ y\sqrt{\tilde D}}2
\in \O_{\tilde F}$  with $x$ and $y$ odd such that
$$
\frac{\tilde\Delta}4 \equiv (\frac{x+y\sqrt{\tilde D}}2)^2 \mod 4
\mathfrak l,
$$
i.e.,
$$
a + 4 D b_1 + 2D^2 c_1 + y^2 \frac{\tilde D -n^2}4
 +(\frac{x -\mu n y}2)^2 + (xy +1) \frac{\mu
n +\sqrt{\tilde D}}2\equiv 0 \mod 4 \mathfrak l.
$$
  Notice that $\frac{\tilde D
-n^2}{4} =D (2ac_1 -4b_1^2)$. So $\tilde K =\tilde
F(\sqrt{\tilde{\Delta}})$ is split at $\mathfrak l$ if and only if
there are odd integers $x$ and $y$ such that
$$
a + 4D(b_1-b_1^2) + 2 D c_1 (D+a) + (\frac{x -\mu n y}2)^2 + (xy +1)
\frac{\mu n +\sqrt{\tilde D}}2\equiv 0 \mod 4 \mathfrak l.
$$
That is
$$
a+ (\frac{x +a y}2)^2 + (xy +1) \frac{\mu n +\sqrt{\tilde D}}2
\equiv 0 \mod 4 \mathfrak l,
$$
since $D +a \equiv 0 \mod 4$ and $-\mu n \equiv a \mod 4$. In
particular, one has to have $x + a y \equiv 2 \mod 4$, and thus $xy
\equiv -1 \mod 4$ (recall that $a \equiv -1 \mod 4$). So the above
congruence is equivalent to
$$
a + (\frac{x +a y}2)^2 \equiv 0 \mod 4 \mathfrak l,
$$
i.e.,  $(-a, 2)_2=1$. So $\tilde K/\tilde F$ is split if and only if
$(-a, 2)_2 =1$.

{\bf Subcase 2}: Now we assume that $c =4 c_1$ with $c_1$ odd. In
this case, $b=2b_1$ with $b_1$ odd, and $\frac{\tilde\Delta}4 \equiv
0 \mod \mathfrak l$,  not easy to deal with. We switch  to
${\tilde\Delta}'$. We have
\begin{equation} \label{eqI4.10}
\frac{{\tilde\Delta}'}4 = \frac{\mu n +\sqrt{\tilde D}}2 - D c_1
\not\equiv 0 \mod \mathfrak l.
\end{equation}
 Since $\tilde\Delta
{\tilde\Delta}' = D v^2$ if $\Delta = \frac{u +v\sqrt D}2$,
$\tilde K=\tilde F(\sqrt{\tilde\Delta}) =\tilde
F(\sqrt{D{\tilde\Delta}'})$ is split at $\mathfrak l$ if and only
if there is odd integers $x$ and $y$ such that
$$
\frac{D {\tilde\Delta}'}4  \equiv  (\frac{x + y \sqrt{\tilde D}}2)^2
\equiv \mod 4 \mathfrak l.
$$
By (\ref{eqI4.9}) and (\ref{eqI4.10}), this is equivalent to
$$
-c_1 \equiv (\frac{x -\mu n y}2)^2 +D(ac_1-b_1^2) + (x y -D)
\frac{\mu n +\sqrt{\tilde D}}2 \mod 4 \mathfrak l.
$$
Since
$$
-a +c_1 + D(ac_1 -b_1^2) =D(a+1)(c_1-1)+ (D-1)(a-c_1) + D(1-b_1^2)
 \equiv 0 \mod 4 \mathfrak l,
 $$
the above congruence is equivalent to
\begin{equation} \label{eqI4.11}
-a \equiv (\frac{x -\mu  n y}2)^2 + (x y -D) \frac{\mu +\sqrt{\tilde
D}}2 \mod 4 \mathfrak l.
\end{equation}
In particular, $\frac{x -\mu  n y}2$ has to be odd. Since
$$
- \mu n = a + 2 D b_1 + (D^2 -D) c_1 \equiv a+2 \mod 4,
$$
this means $\frac{x+ a y}2 \equiv 0 \mod 2$, and thus $xy \equiv 1
\mod 4$. So (\ref{eqI4.11}) is equivalent to
$$
-a \equiv (\frac{x -\mu  n y}2)^2  \mod 4 \mathfrak l.
$$
Therefore, $\tilde K/\tilde  F$ is split at $\mathfrak l$ if and
only if $(-a, 2)_2 =1$. This finally finishes the proof of the
lemma.
\end{proof}

\begin{theo} \label{theoI4.4} One has
\begin{equation} \label{eqI2.7}
b_1(p) = \sum_{0 < n <\sqrt{\tilde D}, \frac{\tilde D -n^2}{4D}
\in p \mathbb Z_{>0}} (\ord_p\frac{\tilde D -n^2}{4D} +1)
\sum_{\mu} \beta(p, \mu n),
\end{equation}
where $\beta(p, \mu n)$ is given as in Theorem \ref{newtheo1.3}.
\end{theo}

\begin{proof} The  discussion between Lemma \ref{lemI4.1}  and Lemma \ref{lemI4.2}  gives
(\ref{eqI2.7}) with
$$
b(p,\mu n) =\prod_{l} b_l(p, \mu n).
 $$
 Here  for $l \ne p$
$$
b_l(p, n) = \prod_{\mathfrak l | l}\rho_{\mathfrak l}(t_{\mu n}
d_{\tilde K/\tilde F} \mathfrak p^{-1}).
$$
For $l=p$, $p\O_{\tilde F} = d_{\tilde K/\tilde F} d_{\tilde
K/\tilde F}'$. Write $\mathfrak p=d_{\tilde K/\tilde F}'$, then
$$
b_p(p,\mu n) =\begin{cases}
  0 &\ff {\mathfrak p} \hbox{ split  in } \tilde K,
  \\
   \rho_{\mathfrak p}(t_{\mu n}
d_{\tilde K/\tilde F} \mathfrak p^{-1})
  &\ff \mathfrak p \hbox{ not
split  in } \tilde K.
\end{cases}
$$
 When $l \nmid \frac{\tilde D -n^2}{4D}$, one has clearly $b_l(p, \mu n)=1$.

When $l=p$,
\begin{align*}
b_p(p, \mu n) &= \begin{cases} 0 &\ff \mathfrak p \hbox { is split
in } \tilde K,
\\
 \frac{1+ (-1)^{t_p -1}}2 &\ff \mathfrak p \hbox { is  not split
in } \tilde K,
\end{cases}
\\
 &=\frac{1 -(-\alpha_p, p)_p^{t_p}}2
 \\
  &=\beta_p(p, \mu n)
 \end{align*}
 as claimed.

  When $l
|\frac{\tilde D -n^2}{4D}$, but $l \ne p$,  $l$ is split in $\tilde
F$ or $ l | \tilde D$ is ramified in $\tilde K$, In the split case,
we choose the splitting $l\O_{\tilde F} =\mathfrak l \mathfrak l'$
so that (\ref{eqI2.5}) holds, and  $\rho_{\mathfrak l'}(t_{ n}
d_{\tilde K/\tilde F} \mathfrak p^{-1})=1$. In either case, $b_l(p,
\mu n) =\rho_{\mathfrak l}(t_{\mu n} d_{\tilde K/\tilde F} \mathfrak
p^{-1})$. Now Lemma \ref{lemI4.2} and (\ref{eqI2.4}) give the
desired formula for $b_l(p, \mu n)$.
\end{proof}

{\bf Proof of Theorem \ref{maintheo}}: Now Theorem \ref{maintheo}
follows from Theorems \ref{newtheo1.3} and \ref{theoI4.4}.

\section{Proof of Theorem \ref{theo1.4}}
\label{newsect7}

A holomorphic Hilbert modular form for $\SL_2(\O_F)$  is called {\em
normalized integral}, if all its Fourier coefficients at the cusp
$\infty$ are rational integers with greatest common divisor $1$.

\begin{lem} \label{lem7.1} Assume $D=5, 13$ or $17$. Then for every integer $m>0$,
there is a positive integer $a(m) >0$ and a normalized integral
holomorphic Hilbert modular form $\Psi_m$ such that
$$
\div \Psi_m =a(m) \mathcal T_m.
$$
\end{lem}
\begin{proof} Let $S_2^+(D, (\frac{D}{}))$ be the space of elliptic
modular forms of weight $2$, level $D$, and Nebentypus character
$(\frac{D}{})$ such that its Fourier coefficients satisfy $a(n)=0$
if $(\frac{D}{n})=-1$. Then a well-known theorem of Hecke asserts
$\dim S_2^+(D, (\frac{D}{})) =0$ for primes  $D=5, 13, 17$. By a
Serre duality theorem of Borcherds \cite{Bo2} and Borcherds's
lifting theorem  \cite{Bo1} (see \cite{BB} in our special setting),
there is Hilbert modular form $\Psi_m$ such that $\div
\Psi_m(\mathbb C) =T_m$ and sufficient large power of $\Psi_m$ is a
normalized integral Hilbert modular form. Replacing $\Psi_m$ by a
sufficient large power if necessary we may assume that $\Psi_m$ is a
normalized integral holomorphic Hilbert modular form. So $\div
\Psi_m$ is flat over $\mathbb Z$ and thus $\div \Psi_m =a(m)
\mathcal T_m$.
\end{proof}

{\bf Proof of Theorem \ref{theo1.4}}: Let $\hat{\omega}=(\omega,
\|\,\|_{\Pet})$ be the metrized Hodge bundle on $\tilde{\mathcal
M}$ with the Petersson metric defined in Section \ref{newsect2}.
Let $\tilde{\mathcal T}_1$ be the closure of $\mathcal T_1$ in
$\tilde{\CalM}$. Let $\Psi_1$ be a normalized integral Hilbert
modular form of weight $c(1)$ given in Lemma \ref{lem7.1}. Then
$\Psi_1$ can be extended to a section of $\omega^{c(1)}$, still
denoted by $\Psi_1$  such that
$$
\div \Psi =a(1) \tilde{\mathcal T}_1.
$$
 Since $\CM(K)$ never intersects with the
boundary $\tilde{\CalM}-\CalM$, $\tilde{\mathcal T}_1.\CM(K) =
\mathcal T_1.\CM(K)$. 
So
\begin{align*}
c(1) \Ht_{\hat{\omega}}(\CM(K))
 &=\Ht_{\widehat{\div}(\Psi_1)} (\CM(K))
 \\
 &= a(1) \CM(K).\mathcal T_1 -\frac{2}{W_K} \sum_{z \in
 \operatorname{CM}(K)} \log ||\Psi_1(z)||_{\Pet}
 \\
  &=\frac{a(1)}{2} b_1 -\frac{a(1)}{2} \frac{W_{\tilde K}}{W_K} b_1
    + \frac{c(1)}2 \frac{W_{\tilde K}}{W_K}\Lambda(0, \chi_{\tilde K/\tilde F}) \beta(\tilde K/\tilde F)
  \end{align*}
by Theorem \ref{maintheo} and \cite[Theorem 1.4]{BY}. It is not hard
to check that
$$
W_K=W_{\tilde K} =\begin{cases}
   10 &\ff K=\tilde K=\mathbb Q(\zeta_5),
   \\
    2 &\hbox{otherwise}.
\end{cases}
$$
Let $M=K\tilde K$ be the Galois closure of $K$ (and $\tilde K$) over
$\mathbb Q$, view both $\chi_{\tilde K/\tilde F}$ and $\chi_{K/F}$
as characters of $\Gal(M/\tilde F)$ and $\Gal(M/F)$ respectively by
class field theory.  Then
$$
\pi =\Ind_{\Gal(M/\tilde F)}^{\Gal(M/\mathbb Q)} \chi_{\tilde
K/\tilde F} =\Ind_{\Gal(M/F)}^{\Gal(M/\mathbb Q)} \chi_{ K/ F}
$$
is the unique two dimensional irreducible representation of
$\Gal(M/\mathbb Q)$ when $K$ is not cyclic (when $K$ is cyclic, the
identity is trivial).  So

$$L(s, \chi_{\tilde K/\tilde F}) = L(s, \chi_{K/F}) =L(s, \pi),
$$
and thus $\beta(\tilde K/\tilde F) =\beta(K/F)$. Finally,
\cite[(9.2)]{BY} asserts
$$
\Lambda(0, \chi_{\tilde K/\tilde F})=\frac{2
\#\mathrm{CM}(K)}{W_K}.
$$
So
$$
h_{\hat{\omega}}(\CM(K)) =\frac{\#\mathrm{CM}(K)}{W_K} \beta(K/F).
$$
Combining this with (\ref{neweq2.7}), one obtains
\begin{equation}
h_{\mathrm{Fal}} (A)=\frac{1}{2} \beta(
 K/F).
 \end{equation}
This proves Theorem \ref{theo1.4}.


\begin{thebibliography}{EMOT}\label{secref}

\bibitem[Ad]{Ad} {\em  G. Anderson},   Logarithmic derivatives of Dirichlet $L$-functions and the periods of abelian varieties.
  Compositio Math.  {\bf 45}(1982), 315--332.


\bibitem[Bo1]{Bo1} {\em R. E. Borcherds}, Automorphic forms with
  singularities on Grassmannians, Invent. Math. {\bf 132} (1998),
  491--562.

\bibitem[Bo2]{Bo2} {\em R. E. Borcherds}, The Gross-Kohnen-Zagier
  theorem in higher dimensions, Duke Math. J. {\bf 97} (1999),
  219--233.


\bibitem[BB]{BB} {\em J. H. Bruinier and M. Bundschuh}, On Borcherds products associated with lattices of
prime discriminant, Ramanujan J. {\bf 7} (2003), 49--61.


\bibitem[BBK]{BBK} {\em J. Bruinier, J. Burgos Gill, and U. K\"uhn},
 Borcherds products and arithmetic intersection theory on Hilbert
 modular surfaces, Duke Math. J. {\bf 139} (2007), 1-88.



\bibitem[BY]{BY} {\em J. H. Bruinier and T. Yang}, CM values of
Hilbert modular functions, Invent. Math. {\bf 163} (2006), 229--288.




 \bibitem[Co]{Col} {\em P. Colmez}, P\'eriods des vari\'et\'es
 ab\'eliennes \`a multiplication complex, Ann. Math., 138(1993),
 625-683.

\bibitem[DP]{DP} {\em P. Deligne and G. Pappas}, Singularit´es des espaces de modules de Hilbert, en les
caract´eristiques divisant le discriminant, Compos. Math. 90
(1994), 59–79.

\bibitem[Fa]{Faltings} {\em G. Faltings}, Finiteness theorems for abelian varieties over number fields.
Translated from the German original [Invent. Math. 73 (1983),
349--366; ibid. 75 (1984),  381; ] by Edward Shipz. Arithmetic
geometry (Storrs, Conn., 1984),  9-27, Springer, New York, 1986.

\bibitem[Go]{Go} {\em E. Goren}, Lectures on Hilbert modular
varieties and modular forms, CRM monograph series 14, 2001.

\bibitem[Gr]{Gr} {\em B. Gross}, On the periods of abelian integrals and a formula of Chowla and Selberg.
With an appendix by David E. Rohrlich.  Invent. Math.  {\bf 45}
(1978), 193-211.

\bibitem[GK]{GK}  {\em B. Gross and K. Keating}, On the intersection
of modular correspondences, Invent. Math. {\bf 112}(1993), 225-245.

\bibitem[GZ1]{GZ2} {\em B. Gross and D. Zagier}, On singular moduli.  J. Reine Angew. Math.  {\bf 355 } (1985), 191-220.

\bibitem[GZ2]{GZ} {\em B. Gross and D. Zagier}, Heegner points and derivatives of $L$-series.  Invent. Math. {\bf 84}  (1986),  225-320.


\bibitem[HZ]{HZ} {\em F. Hirzebruch and D. Zagier}, Intersection Numbers of Curves on Hilbert Modular
Surfaces and Modular Forms of Nebentypus, Invent. Math. {\bf 36}
(1976), 57--113.


\bibitem[Ig]{Ig} {\em J.-I. Igusa}, Arithmetic Variety of Moduli for Genus Two.
Ann. Math. 72, (1960), 612-649



\bibitem[La]{La} {\em K. Lauter},
Primes in the denominators of Igusa Class Polynomials, preprint,
pp3, http://www.arxiv.org/math.NT/0301240/.

\bibitem[Ku1]{KuAnnal} {\em S. Kudla}, Central derivatives of Eisenstein series and height pairings.
 Ann. of Math. (2)  {\bf 146 } (1997),  545-646.


\bibitem[Ku2]{Ku2} {\em S. Kudla},   Special cycles and derivatives of Eisenstein series, in  Heegner points and Rankin L-series,
243-270, Math. Sci. Res. Inst. Publ., 49, Cambridge Univ. Press,
Cambridge, 2004.

\bibitem[KR1]{KR1} {\em S. Kudla and M. Rapoport},   Arithmetic Hirzebruch Zagier cycles, J. reine angew. Math., {\bf 515 }(1999),
155-244.

\bibitem[KR2]{KR2} {\em S. Kudla and M. Rapoport}, Cycles on Siegel 3-folds and derivatives of Eisenstein series,
 Annales Ecole. Norm. Sup. {\bf 33} (2000), 695-756.

\bibitem[KRY1]{KRY1}  {\em S. Kudla, M. Rapoport, and T.H. Yang},
 Derivatives of Eisenstein Series and Faltings heights,   Comp. Math., {\bf 140 }(2004),   887-951.

\bibitem[KRY2]{KRY2}  {\em S. Kudla, M. Rapoport, and T.H. Yang},  Modular forms and special cycles on
Shimura curves,  Annals of Math. Studies series, vol 161, 2006,
Princeton Univ. Publ.


\bibitem[Se]{Se} {\em J.-P. Serre}, A course in Arithmetic, GTM
{\bf 7}, Springer-Verlag, New York, 1973.

\bibitem[Vo]{Vo} {\em I. Vollaard},  On the Hilbert-Blumenthal moduli
problem,
  J. Inst. Math. Jussieu  4  (2005),
653--683.

\bibitem[We1]{We2} {\em T. Wedhorn}, Caculation of representation densities, Chapter 15
in ARGOS seminar on Intersections of Modular Correspondences, p.
185-196, to appear in Asterisque.

\bibitem[We2]{We1} {\em T. Wedhorn}, The genus of the endomorphisms
of a supersingular elliptic curve, Chapter 5 in ARGOS seminar on
Intersections of Modular Correspondences, p. 37-58, to appear in
Asterisque.



\bibitem[Ya1]{YaDensity1} {\em T.H. Yang},
  An explicit formula for local densities of quadratic
  forms, J. Number Theory  {\bf 72}(1998), 309-356.

\bibitem[Ya2]{YaDensity2} {\em T.H. Yang}, Local densities of $2$-adic quadratic
forms,  J. Number Theory  {\bf 108}(2004), 287-345.

\bibitem[Ya3]{Ya3} {\em T. H. Yang}, Chowla-Selberg Formula and Colmez's Conjecture, preprint, 2007, pp17.


\bibitem[Ya4]{Ya4} {\em  T. H. Yang}, Arithmetic Intersection on a Hilbert Modular Surface and  Faltings' Height,
preprint, 2008,  pp45.

\bibitem[Yu]{Yu} {\em C.F. Yu},  The isomorphism classes of abelian
varieties of  CM types, Jour. pure and Appl. Algebra, {\bf
187}(2004), 305-319.

\bibitem[Zh1]{Zh1} {\em S.W. Zhang},
Heights of Heegner cycles and derivatives of L-series, Invent. Math.
{\bf 130} (1997),  99--152.

\bibitem[Zh2]{Zh2} {\em S. W. Zhang},
Heights of Heegner points on Shimura curves,
 Ann. of Math. (2) {\bf 153} (2001),  27--147.

\bibitem[Zh3]{Zh3} {\em S. W. Zhang},
Gross-Zagier formula for GL(2), Asian. J. Math. 5 (2001), 183--290;
II, Heegner points and Rankin L-series, 191--214, Math. Sci. Res.
Inst. Publ., 49, Cambridge Univ. Press, Cambridge, 2004.



\end{thebibliography}
\end{document}